\documentclass[letterpaper,    
               english,        
               11pt,           
               twoside]        
               {amsart}        
\usepackage[english]{babel}                   
\usepackage{color}
\usepackage[letterpaper,                      
            bookmarks,                        
            bookmarksopen=true,               
            bookmarksnumbered=true,           
            pdfauthor={Christof Puhle},       
            pdftitle={The Killing spinor      
                      equation with higher
                      order potentials},
            pdfkeywords={Connections with torsion,
                         Killing spinor equation, G-structures,
                         type II string theory},   
            colorlinks,                       
            linkcolor=black,                  
            citecolor=black,                  
            urlcolor=black,                   
            plainpages=false]                 
            {hyperref}                        
\usepackage[nobysame]
           {amsrefs}                          
\usepackage{amstext}                          
\usepackage{amssymb}                          
\usepackage{pifont}                           
\advance\oddsidemargin by -0.93cm
\advance\evensidemargin by -0.92cm
\advance\topmargin by -0.55cm
\theoremstyle{plain}

\newtheorem{lem}{Lemma}[section]
\newtheorem{thm}{Theorem}[section]            
\newtheorem{prop}{Proposition}[section]
\theoremstyle{definition}
\newtheorem{exa}{Example}[section]
\newtheorem{rmk}{Remark}[section]

\newtheorem*{thank}{Thanks}
\newcommand{\bdm}{\begin{displaymath}}
\newcommand{\edm}{\end{displaymath}}
\newcommand{\be}[1][*]{\begin{equation#1}}
\newcommand{\ee}[1][*]{\end{equation#1}}
\newcommand{\bea}[1][*]{\begin{eqnarray#1}}
\newcommand{\eea}[1][*]{\end{eqnarray#1}}
\newcommand{\ba}[1]{\begin{array}{#1}}
\newcommand{\ea}{\end{array}}
\newcommand{\btab}{\begin{tabular}}
\newcommand{\etab}{\end{tabular}}
\newcommand{\bcen}{\begin{center}}
\newcommand{\ecen}{\end{center}}
\def\haken{\mathbin{
                 \hskip-0.2pt                
                 \vrule height0.4pt width5.0pt depth0pt   
                 \kern-0.4pt
                 \vrule height6.0pt width0.4pt depth0pt
                 \hskip0.2pt  }}
\def\shaken{\mathbin{
                 \hskip-0.2pt                  
                 \vrule height0.3pt width3.8pt depth0pt   
                 \kern-0.3pt
                 \vrule height4.7pt width0.3pt depth0pt
                 \hskip0.2pt }}                         
\newcommand{\x}{\times}
\newcommand{\op}{\oplus}
\newcommand{\ox}{\otimes}

\newcommand{\ra}{\rightarrow}

\newcommand{\qqs}{\forall}

\newcommand{\Id}{\ensuremath{\mathrm{Id}}}

\newcommand{\C}{\ensuremath{\mathbb{C}}}
\newcommand{\R}{\ensuremath{\mathbb{R}}}

\newcommand{\M}{\ensuremath{\mathcal{M}}}

\newcommand{\W}{\ensuremath{\mathcal{W}}}
\newcommand{\xhi}{\ensuremath{\mathcal{X}}}

\newcommand{\G}{\ensuremath{\mathrm{G}}}

\newcommand{\mc}[1]{\mathcal{#1}}
\newcommand{\mrm}[1]{\mathrm{#1}}
\newcommand{\vphi}{\ensuremath{\varphi}}    

\let\mathg\mathfrak

\newcommand{\hut}{\wedge}

\newcommand{\Ric}{\ensuremath{\mathrm{Ric}}}
\newcommand{\Scal}{\ensuremath{\mathrm{Scal}}}

\newcommand{\Hol}{\ensuremath{\mathrm{Hol}}}
\newcommand{\hol}{\ensuremath{\mathg{hol}}}
\newcommand{\Iso}{\ensuremath{\mathrm{Iso}}}

\newcommand{\diag}{\ensuremath{\mathrm{diag}}}

\newcommand{\SL}{\ensuremath{\mathrm{SL}}}
\newcommand{\un}{\ensuremath{\mathg{u}}}
\newcommand{\su}{\ensuremath{\mathg{su}}}
\newcommand{\SU}{\ensuremath{\mathrm{SU}}}
\newcommand{\U}{\ensuremath{\mathrm{U}}}

\newcommand{\so}{\ensuremath{\mathg{so}}}
\newcommand{\SO}{\ensuremath{\mathrm{SO}}}

\newcommand{\g}{\ensuremath{\mathfrak{g}}}

\newcommand{\m}{\ensuremath{\mathfrak{m}}}

\begin{document}
\thispagestyle{empty}
\date{January 14, 2008}
\title{The Killing spinor equation with higher order potentials}
\author{Christof Puhle}
\address{\hspace{-5mm} 
{\normalfont\ttfamily puhle@mathematik.hu-berlin.de}\newline
Institut f\"ur Mathematik \newline
Humboldt-Universit\"at zu Berlin\newline
Unter den Linden 6\newline
10099 Berlin, Germany}
\thanks{Supported by the SFB 647: `Space--Time--Matter' and the SPP 1154:
`Global Differential Geometry'}
\subjclass[2000]{Primary 53 C 25; Secondary 81 T 30}
\keywords{Connections with torsion, Killing spinor equation, type II string theory}
\begin{abstract}
Let $(M^n,g)$ be a Riemannian spin manifold. The basic equations in supergravity
models of type IIa string theory with $4$-form flux involve a $3$-form $T$, a
$4$-form $F$, a spinorial covariant derivative $\nabla$ depending on $\nabla^g$,
$T$, $F$, and a $\nabla$-parallel spinor field $\Psi$. We classify and construct
many explicit families of solutions to this system of spinorial field equations
by means of non-integrable special geometries. The latter include $\alpha$-Sasakian
structures in dimensions $5$ and $7$, almost Hermitian structures in dimension $6$
and cocalibrated $\G_2$-structures in dimension $7$. We show that there are several
examples also satisfying an additional constraint for the energy-momentum tensor.
\end{abstract}
\maketitle
\setcounter{tocdepth}{1}
\begin{center}
\begin{minipage}{0.7\linewidth}
    \begin{small}
      \tableofcontents
    \end{small}
\end{minipage}
\end{center}
\pagestyle{headings}
%
%
%
%
\section{Introduction}\noindent
In the seventies A.~Gray considered non-integrable special Riemannian
geometries of small dimensions ($n\leq8$). A decade later these structures
started to play a crucial role in the investigations of T.~Friedrich et
al.~concerning eigenvalue estimates for the Dirac operator on a Riemannian
manifold.

Following a latent period, the interest in non-integrable geometries
emerged once again from the developments of theoretical physics in the
common sector of type II string theory. At the turn of the century
T.~Friedrich and his collaborators developed a new and systematic
approach to non-integrable special geometries which takes up certain
aspects of string theory straightforwardly. This approach relies in a
natural way on the notion of characteristic connection $\nabla^c$. This
is an affine connection with totally skew-symmetric torsion $T^c$,
the so-called \emph{torsion form}, which can be associated to certain
$\G$-structures. This point of view lead to show that structures with
parallel torsion ($\nabla^cT^c=0$) are of particular interest. They
provide in fact a basis to solve Strominger's equations in a natural
manner (see \cites{Str86,FI02}).

Since a few years now a more general system of spinorial field equations
than Strominger's has become central in type II string theory. These models
of supergravity -- the so-called models with fluxes -- can be described
geometrically by a tuple $(M^n,g,T,F,\Psi)$ consisting of a Riemannian spin
manifold $(M^n,g)$, a $3$-form $T$, a $4$-form $F$ and a spinor field $\Psi$
satisfying
\be
(\aleph)\quad\quad
\begin{gathered}
\nabla^{g}_X \Psi + \frac{1}{4}\,(X\haken T)\cdot\Psi 
    + p\,(X\haken F)\cdot\Psi + q\,(X\hut F)\cdot\Psi = 0, \\
\Ric^g_{ij} - \frac{1}{4}\,T_{imn}T_{jmn}=0,\quad \delta(T)=0
\end{gathered}
\ee
where $\nabla^g$ denotes the Levi-Civita connection of $(M^n,g)$ and $p,q\in\R$
are real parameters (see article \cite{Duf02}). The first of these three equations
-- the so-called \emph{Killing spinor equation} -- should be satisfied for any
vector field $X\in TM$. If one introduces the new spinorial covariant derivative
\be
\nabla_X\Psi:=\nabla^{g}_X \Psi + \frac{1}{4}\,(X\haken T)\cdot\Psi 
                 + p\,(X\haken F)\cdot\Psi + q\,(X\hut F)\cdot\Psi,
\ee
the Killing spinor equation takes the particularly simple form $\nabla\Psi=0$.
Considering the Kaluza-Klein reduction of $\mc{M}$-theory (see
\cites{WNW85,Ali01,BDS01}) the relevant dimension for $M^n$ lies
between $4$ and $8$. Moreover, additional algebraic constraints occur,
for example the algebraic intertwining between the $3$-form or the
$4$-form and the spinor field $\Psi$
\bea
(\aleph\aleph)\quad\quad&T\cdot\Psi = \lambda\cdot\Psi,\quad
F \cdot \Psi = \kappa\cdot\Psi, \quad \lambda,\kappa\in\C.
\eea
Obviously, the system consisting of $(\aleph)$ and $(\aleph\aleph)$
generalizes Strominger's model by introducing the new degree of
freedom given by a 4-form $F$, usually called \emph{flux form}.\medskip\\\noindent
The task of the present work is the construction of solutions to
the system $(\aleph)$ \& $(\aleph\aleph)$.\medskip\\\noindent
The paper is structured as follows: In \autoref{sec:2} we fit the problem
into the framework of special geometries. We then specialize to the structures
of concern in \autoref{sec:3}. In \autoref{sec:4} we present the classification
and construction techniques used to solve the Killing spinor equation and we
list results thus obtained in sections \ref{sec:5}--\ref{sec:7}. Eventually we
review these in the light of the entire system $(\aleph)$ \& $(\aleph\aleph)$.
%
%
%
%
\section{Setup}\label{sec:2}\noindent
Two questions underlie the whole discussion:
\begin{itemize}
\item Which is the `correct way' to obtain solutions to the entire
      system?\vspace{1mm}
\item Can the differential forms $T$ and $F$ be chosen or fixed in
      a canonical way?
\end{itemize}
The starting point for answering the first question is to consider the
Killing spinor equation, $\nabla\Psi=0$. Solving this is the central task of
the present paper and will be addressed in the following sections in a systematic
way.

Behind the second question hides the aim to control the high degree of freedom
in the choice of differential forms. To tackle this we study certain classes
of non-integrable $\G$-structures with characteristic connection $\nabla^c$ and
parallel torsion form, $\nabla^cT^c=0$. These structures have been investigated
a lot over the last years (see \cites{FI02,FI03a,Fri06,Sch07}). This approach
has the advantage that many geometric properties follow from the parallelism of
$T^c$, and a natural ansatz for the $3$-form is to require $T$ to be fixed up
to a real parameter: $T\sim T^c$. Unfortunately, there is no natural/canonical
$4$-form in this setup, to the effect that -- at least in principle -- $F$ is
completely arbitrary. To overcome this problem we make a special assumption on
$F$, and furthermore demand it to be parallel with respect to $\nabla^c$.

Precisely, we take a class $(M^n,g,\nabla^cT^c=0)$ of non-integrable $\G$-structures
with parallel torsion, fix a spin structure, define $\nabla^c$-parallel $4$-forms
$F_i$ and assume the following for the differential forms $T$ and $F$:
\be
F=\sum_i A_i\cdot F_i,\quad T=B\cdot T^c, \quad A_i,B \in \R.
\ee
Given this, we try to solve $\nabla\Psi=0$ on the underlying special
structure. Which are the classes of special geometries dealt with
is readily explained.
%
%
%
%
\section{Special geometries}\label{sec:3}\noindent
Our treatise includes the dimensions $5$, $6$ and $7$. In dimension $5$
we investigate $\alpha$-Sasakian structures. As far as dimension $6$ is
concerned we choose almost Hermitian structures with parallel torsion,
classified by N.~Schoemann (see \cites{Sch06,Sch07}). In dimension $7$ we
consider both $\alpha$-Sasakian structures and cocalibrated $\G_2$-structures
with parallel torsion, the latter exhaustively described in \cite{Fri06}.
\subsection{\texorpdfstring{$\alpha$}{Alpha}-Sasakian structures}
We begin with some basic definitions of contact geometry. The book of Blair
\cites{Bla76} and the articles \cites{AG86,CG90} may serve as general references.
An \emph{almost contact metric structure} consists of an odd-dimensional
manifold $M^{2k+1}$ equipped with a Riemannian metric $g$, a vector field $\xi$
of length one, its dual 1-form $\eta$ (\emph{contact form}) as well as an
endomorphism $\vphi$ of the tangent bundle such that the algebraic relations
\be
\vphi(\xi) = 0, \qquad \vphi^2 = -\Id+\eta\ox\xi, \qquad 
g\left( \vphi(X),\vphi(Y)\right) = g(X,Y)-\eta(X)\,\eta(Y)
\ee
are satisfied. The \emph{fundamental form} $\Phi$ and the \emph{Nijenhuis tensor} $N$
of an almost contact metric structure are defined by
\be
\Phi(X,Y) := g\left( X,\vphi(Y)\right), \qquad
N(X,Y) := [\vphi,\vphi](X,Y)+d\eta(X,Y)\cdot\xi.
\ee
There are many special types of almost contact metric structures in the literature.
We introduce those appearing in this paper. An almost contact metric structure
is called \emph{normal} if its Nijenhuis tensor vanishes, $N=0$. A \emph{quasi-Sasakian}
structure has additionally closed fundamental form, $d\Phi=0$. A normal almost
contact metric structure with
\be
N=0,\quad d\eta=\alpha\cdot \Phi, \quad \alpha\in\R\backslash\{0\}
\ee
is called $\alpha$-\emph{Sasakian}. Taking $\alpha=2$ above restricts to
\emph{Sasakian} structures.

It is known (cf.~\cite{FI02}*{Thm.~8.2}) that every quasi-Sasakian manifold
$(M^{2k+1},g,\xi,\eta,\vphi)$ admits a unique metric connection $\nabla^c$
with totally skew-symmetric torsion $T^c$,
\be
g\left( \nabla^{c}_{X} Y, Z \right)  =   g\left( \nabla^{g}_{X} Y, Z \right) 
                                        + \frac{1}{2}\cdot T^c\left(X,Y,Z\right),
\ee
preserving the quasi-Sasakian structure, $\nabla^c\xi=\nabla^c\vphi=0$. The
torsion form is given by
\be
T^c = \eta\hut d\eta.
\ee
Obviously, $T^c$ is parallel with respect to $\nabla^c$ in the case
of $\alpha$-Sasakian structures.

For every almost contact metric structure $(M^{2k+1},g,\xi,\eta,\vphi)$
there exists an oriented orthonormal frame $(e_1,\ldots,e_{2k+1})$ such
that
\be
\xi = e_{2k+1}, \quad \Phi = e_{12} + e_{34} + \ldots + e_{2k-1,2k}.
\ee
We will call this frame an \emph{adapted frame} of the corresponding structure.
Here and henceforth we shall not distinguish between vectors and covectors
and use the notation $e_{i_1\ldots i_m}$ for the exterior product $e_{i_1}\hut
\ldots\hut e_{i_m}$.
\subsection{Almost Hermitian structures}
An \emph{almost Hermitian structure} $(M^6,g,J)$ is a six-dimen\-sional manifold $M^6$
equipped with a Riemannian metric $g$ and an orthogonal almost complex structure
$J: TM\ra TM$,
\be
J^2 = -\Id_{TM} , \quad g(JX,JY) = g(X,Y).
\ee
The \emph{K\"ahler form} $\Omega$ is defined by $\Omega (X,Y) := g(JX,Y)$. A
\emph{nearly K\"ahler} structure is an almost Hermitian structure that satisfies
the condition
\be
(\nabla^g_XJ)(X)=0, \quad \qqs X\in TM.
\ee
Finally, \emph{K\"ahler} manifolds are characterized by the $\nabla^g$-parallelism
of the almost complex structure, $\nabla^gJ=0$.

A six-dimensional almost Hermitian structure can be understood as a $\U(3)$-reduction
of the corresponding bundle of orthonormal frames $\mc{F}(M^{6},g)$. The space of
$3$-forms
\be
\Lambda^3(\R^6)=\Lambda^3_2(\R^6)\op\Lambda^3_{12}(\R^6)\op\Lambda^3_6(\R^6)
\ee
decomposes into three irreducible $\U(3)$-components as described for example in
\cite{AFS05}. Decomposing the Lie algebra $\so(6)$ into $\un(3)$ and its orthogonal
complement $\m^6$ the fundamental classes of six-dimensional almost Hermitian
structures are defined by the irreducible $\U(3)$-sub\-modules of
\be
\R^6\ox\m^6=\W_1\op\W_2\op\W_3\op\W_4.
\ee
$\W_1$, $\W_3$ and $\W_4$ can be characterized in terms of the irreducible
submodules of $\Lambda^3(\R^6)$ above (see \cite{AFS05}). The notation is
in accordance with the Gray-Hervella classification \cite{GH80}.

It is known (cf.~\cite{AFS05}*{Cor.~3.5}) that there exists an unique affine
connection $\nabla^c$ with totally skew-symmetric torsion $T^c$ preserving
the almost Hermitian structure ($\nabla^cJ=0$) if and only if this structure
is of type $\W_1\op\W_3\op\W_4$.

We say that six-dimensional almost Hermitian structures of type $\W_1\op\W_3\op
\W_4$ \emph{belong to the class $\mc{C}[\G]$} for the proper subgroup $\G\subset
\U(3)$, if the corresponding torsion form is parallel with respect to the
characteristic connection, $\nabla^cT^c=0$, and if
\be
\Hol(\nabla^c)\subset\G\subset\Iso(T^c).
\ee
Here $\Hol(\nabla^c)$ is the holonomy of the characteristic connection and
$\Iso(T^c)$ the connected component at the identity of the isotropy group
of $T^c$. The classes $\mc{C}[\G]$ were investigated in \cite{Sch07}. We will
discuss those with non-abelian $\G$ in more detail in \autoref{sec:6}.

Just as before, there exists always an \emph{adapted frame} $(e_1,\ldots,e_6)$,
that is such that
\be
\Omega = e_{12} + e_{34} + e_{56}.
\ee
\subsection{Cocalibrated \texorpdfstring{$\G_2$}{G2}-structures}
Consider the space $\R^7$, fix an orientation and denote a chosen
oriented orthonormal basis by $(e_1,\ldots,e_7)$. The Lie group
$\G_2$ can be described as the isotropy group of the $3$-form
\be
(\ast)\quad\quad\omega^3=e_{127}+e_{135}-e_{146}-e_{236}-e_{245}+e_{347}+e_{567}.
\ee
A \emph{$\G_2$-structure} is a triple $(M^7,g,\omega^3)$ consisting
of a seven-dimensional Riemannian manifold $(M^7,g)$ and a $3$-form
$\omega^3$ such that there exists an oriented orthonormal \emph{adapted
frame} $(e_1,\ldots,e_7)$ realizing $(\ast)$ at every point.

Obviously, a $\G_2$-structure is a reduction of the structure group
of orthonormal frames of the tangent bundle to $\G_2$. The space of
$3$-forms decomposes into three irreducible $\G_2$-components (see
for example \cite{FKMS97}),
\be
\Lambda^3(\R^7)=\Lambda^3_1(\R^7)\op\Lambda^3_{27}(\R^7)\op\Lambda^3_7(\R^7).
\ee
Using these components one can describe the irreducible $\G_2$-submodules
$\xhi_1$, $\xhi_3$ and $\xhi_4$ of the decomposition
\be
\R^7\ox\m^7=\xhi_1\op\xhi_2\op\xhi_3\op\xhi_4
\ee
where $\m^7$ is the orthogonal complement of $\g_2$ inside $\so(7)$ (see
\cite{FI02}). In this way we obtain the Fern\'andez-Gray description
\cite{FG82} of non-integrable $\G_2$-structures by differential equations.
For example, a $\G_2$-structure is of type $\xhi_1$, i.e.~a \emph{nearly parallel}
structure, if and only if there exists a real number $\lambda\neq0$ such that
\be
d\omega^3=-\lambda\cdot\ast\omega^3.
\ee
$\G_2$-structures of type $\xhi_1\op\xhi_3$ -- the so-called \emph{cocalibrated}
structures -- are characterized by a coclosed $3$-form,
\be
\delta\omega^3=0.
\ee

The following is known (cf.~\cite{FI02}*{Thm.~4.7}) on the existence of
characteristic connection: There exists an unique affine connection
$\nabla^c$ with totally skew-symmetric torsion $T^c$ preserving the
$\G_2$-structure ($\nabla^c\omega^3=0$) if and only if this structure
is of type $\xhi_1\op\xhi_3\op\xhi_4$.

Cocalibrated, non-nearly parallel $\G_2$-structures with parallel characteristic
torsion,
\be
\nabla^cT^c=0,
\ee
for which the holonomy algebra $\hol(\nabla^c)$ of the characteristic connection
is a proper subalgebra $\g\subset\g_2$, are said to \emph{belong to the class
$\mc{C}[\g]$}. Structures of class $\mc{C}[\g]$ with non-abelian $\g$ were investigated
in \cite{Fri06}. We will return to this in \autoref{sec:7}.
%
%
%
%
\section{Techniques}\label{sec:4}\noindent
To solve $\nabla\Psi=0$ as explained we proceed in two steps:
\begin{enumerate}
\item We `classify' solutions, i.e.~determine necessary conditions for
      $\nabla\Psi=0$ to hold.\vspace{1mm}
\item We construct solutions with the help of (1).

\end{enumerate}
\subsection{Classification Technique}
Fix a Riemannian spin manifold $(M^n,g)$ as in section \ref{sec:2}, denote its
spinor bundle by $\Sigma$ and assume there exists a solution $\Psi\in\Gamma(\Sigma)$
of $\nabla\Psi=0$. The spinor field $\Psi$ is an element in the kernel of the
spinorial curvature tensor
\be
\mrm{R}^{\nabla}(X,Y) := \nabla_X\nabla_Y - \nabla_Y\nabla_X - \nabla_{[X,Y]}.
\ee
By understanding this differential operator as an endomorphism of spinors,
we may define its \textit{first and second contraction $K^{\nabla}(X)$ and
$K^{\nabla}$} via Clifford multiplication relative to an adapted frame
$(e_1\,\ldots,e_n)$,
\be
K^{\nabla}(X) := \sum_{k}e_k\cdot\mrm{R}^{\nabla}(e_k,X), \quad
K^{\nabla} := \sum_{k}e_k\cdot K^{\nabla}(e_k).
\ee
We denote the curvature terms related to $\nabla=\nabla^c$ by $\mrm{R}^c(X,Y)$,
$K^c(X)$ and $K^c$. If we compare $\mrm{R}^{\nabla}(X,Y)$ to $\mrm{R}^c(X,Y)$,
we obtain
\be
\mrm{R}^{\nabla}(X,Y)=\mrm{R}^{c}(X,Y)+\sum\limits_{j=1}^{12}R_j(X,Y)
\ee
where the algebraic terms $R_j(X,Y)$ are given by ($s:=(B-1)/4$)
\bea
R_1(X,Y)&=&s\left(T^c(X,Y)\haken T^c\right),\\
R_2(X,Y)&=&p\left(T^c(X,Y)\haken F\right),\\
R_3(X,Y)&=&q\left(T^c(X,Y)\hut F\right),\\
R_4(X,Y)&=&s^2\left( (X\haken T^c)\cdot(Y\haken T^c)-(Y\haken T^c)\cdot(X\haken T^c) \right),\\
R_5(X,Y)&=&s\, p\left( (X\haken T^c)\cdot(Y\haken F)-(Y\haken F)\cdot(X\haken T^c) \right),\\
R_6(X,Y)&=&s\, p\left( (X\haken F)\cdot(Y\haken T^c)-(Y\haken T^c)\cdot(X\haken F) \right),\\
R_7(X,Y)&=&s\, q\left( (X\haken T^c)\cdot(Y\hut F)-(Y\hut F)\cdot(X\haken T^c) \right),\\
R_8(X,Y)&=&s\, q\left( (X\hut F)\cdot(Y\haken T^c)-(Y\haken T^c)\cdot(X\hut F) \right),\\
R_9(X,Y)&=&p^2\left( (X\haken F)\cdot(Y\haken F)-(Y\haken F)\cdot(X\haken F) \right),\\
R_{10}(X,Y)&=&p\, q\left( (X\haken F)\cdot(Y\hut F)-(Y\hut F)\cdot(X\haken F) \right),\\
R_{11}(X,Y)&=&p\, q\left( (X\hut F)\cdot(Y\haken F)-(Y\haken F)\cdot(X\hut F) \right),\\
R_{12}(X,Y)&=&q^2\left( (X\hut F)\cdot(Y\hut F)-(Y\hut F)\cdot(X\hut F) \right).
\eea
Let
\be
M_j(X) := \sum_{k}e_k\cdot R_j(e_k,X), \quad
M_j := \sum_{k}e_k\cdot M_j(e_k)
\ee
denote the contractions. We derive the following formulae for the first and
second contraction of $\mrm{R}^\nabla(X,Y)$ using the proof of \cite{FI02}*{Cor.~3.2}:
\be
K^{\nabla}(X) = \frac{1}{2}\cdot\Ric^c(X)-\frac{1}{2}\cdot(X\haken \sigma^{T^c})
+\sum\limits_{j=1}^{12}M_j(X),\quad
K^{\nabla} = -\frac{1}{2}\cdot\Scal^c-2\cdot\sigma^{T^c}+\sum\limits_{j=1}^{12}M_j.
\ee
Here $\sigma^{T^c}$ denotes the $4$-form
\be
\sigma^{T^c}:=\frac{1}{2}\cdot\sum_i(e_i\haken T^c)\hut(e_i\haken T^c).
\ee

To summarize: If $T$ and $F$ are chosen parallel with respect to $\nabla^c$,
we are able to compute the first and second contraction of the spinorial
curvature tensor $\mrm{R}^\nabla(X,Y)$ algebraically, provided we know the
characteristic Ricci tensor $\Ric^c$. By a careful inspection of the kernels
of these spinorial endomorphism we shall obtain necessary conditions for
$\nabla\Psi=0$. We will demonstrate this technique and exhibit the corresponding
contractions $K^{\nabla}(X)$ and $K^{\nabla}$ in the case of almost Hermitian
structures of class $\mc{C}[\SU(3)]$. 

For simplicity we split the problem $\nabla\Psi=0$ in the following sections
by considering the covariant derivatives
\be
\nabla=\left\{\begin{array}{ccl}
	 \nabla^0 & \mrm{for} &p=(n-4)/4,\,\,q=1\\
	 \nabla^1 & \mrm{for} &p=(n-4)/4,\,\,q\,\,\mrm{arbitrary}\\
	 \nabla^2 & \mrm{for} &p=0,\,\,q=1
\end{array} \right.
\ee
thus
\be
\nabla^{0}_{X}\Psi=\nabla^{g}_X \Psi + \frac{1}{4}\,(X\haken T)\cdot\Psi 
                 + \frac{n-4}{4}\,(X\haken F)\cdot\Psi + (X\hut F)\cdot\Psi,
\ee
\be
\nabla^{1}_{X}\Psi=\nabla^{g}_X \Psi + \frac{1}{4}\,(X\haken T)\cdot\Psi 
                 + \frac{n-4}{4}\,(X\haken F)\cdot\Psi + q\,(X\hut F)\cdot\Psi,
\ee
\be
\nabla^{2}_{X}\Psi=\nabla^{g}_X \Psi + \frac{1}{4}\,(X\haken T)\cdot\Psi
+(X\hut F)\cdot\Psi.
\ee
Up to a rescaling of $F$ the first case refers to a particular ratio between the
parameters $p$ and $q$ in $\nabla$, namely $4\, p=(n-4)\, q$. This is special in
several ways, so we will call the corresponding equation \emph{of special type}.
For example, this is justified by the fact that the Dirac operator defined by
$\nabla$ does not depend on $F$ for $p/q=(n-4)/4$ (see \cite{AF03}).
\subsection{Construction Technique}
The construction of solutions depends very specifically on the classification
results, and therefore on the underlying special structure. However, there are
two general ideas on how one can use these results for constructing solutions.

Firstly, we restrict further by only considering simply connected spin manifolds,
and suppose that $K^{\nabla}(X)\Psi=0$ is satisfied for every spinor field belonging
to a one-dimensional spin subbundle $\Sigma^1$ of $\Sigma$. The methods of \cite{FK89}
show that if $\mrm{R}^{\nabla}(X,Y)|_{\Sigma^1}\equiv0$ holds and if $\nabla$ preserves
$\Sigma^1$, then there exists a $\nabla$-parallel spinor field in $\Sigma^1$.

Secondly, solving $K^{\nabla}(X)\Psi=0$ suggests in many cases to look for a `special'
spinor field, for example one parallel with respect to $\nabla^c$ or a Killing spinor.
If we can prove the existence of this spinor field, we have a starting point for
the construction of solutions.
\medskip\\\noindent
The presentation of results thus obtained is divided into broad parts for each class of
special geometries introduced in \autoref{sec:3}. We describe examples relative to the
various subclasses and state the classification results. Eventually we discuss which
necessary conditions might be sufficient for the construction of solutions. An exception
is made for almost Hermitian structures of class $\mc{C}[\SU(3)]$. For that case we will
additionally present the techniques for clarity. This also means that we omit the
completely analogous but lengthy proofs for all other cases, to be found in the thesis
\cite{Puh07}.
%
%
%
%
\section{\texorpdfstring{$\alpha$}{Alpha}-Sasakian structures}\label{sec:5}\noindent
Relative to an adapted frame the torsion form $T^c$ of an $\alpha$-Sasakian
structure is given by (see \autoref{sec:3})
\be
T^c=\alpha\cdot\Phi\hut\eta=\alpha\cdot( e_{12} + e_{34} +
\ldots + e_{2k-1,2k} )\hut e_{2k+1}.
\ee
Recall \cite{FI02} that the characteristic Ricci tensor is symmetric,
\be
\Ric^c(X,Y) = \Ric^c(Y,X),
\ee
due to the $\nabla^c$-parallelism of $T^c$.
\begin{lem}
On an $\alpha$-Sasakian manifold $(M^{2k+1},g,\xi,\eta,\vphi)$ of dimension
$2k+1$ the Riemannian Ricci curvature in the direction $\xi$ is equal to $k\,\alpha^2/2$.
\end{lem}
\begin{proof}
In analogy to \cite{Bla76}, the relation
\be
\alpha\, g(X,\vphi(Y)) = d\eta(X,Y) = 2\, g(\nabla^{g}_X\xi,Y)
\ee
shows that $\nabla^g$ acts on $\xi$ by the endomorphism $\vphi$,
\be
\nabla^{g}_X\xi = -\frac{\alpha}{2}\cdot\vphi(X).
\ee
Let $X$ now denote an unit vector field orthogonal to $\xi$. Then
\bea
\mrm{R}^{g}(\xi,X)\xi
& = & -\frac{\alpha}{2}\cdot\nabla^{g}_\xi\vphi(X)+\frac{\alpha}{2}\cdot\vphi([\xi,X])
 =  -\frac{\alpha}{2}\cdot\nabla^{g}_\xi\vphi(X)+\frac{\alpha}{2}\cdot[\xi,\vphi(X)]\\
& = & -\frac{\alpha}{2}\cdot\nabla^{g}_{\vphi(X)}\xi
\, = \, \frac{\alpha^2}{4}\cdot\vphi^2(X)\, = \, -\frac{\alpha^2}{4}\cdot X ,
\eea
and hence $g(\mrm{R}^{g}(\xi,X)X,\xi)=\alpha^2/4$.
\end{proof}\noindent
To compare $\Ric^c$ and $\Ric^{g}$ we can use the following formula (see \cite{FI02}):
\be
\Ric^{c}_{ij}  =  \Ric^{g}_{ij}-\frac{1}{4}\cdot T^{c}_{imn}T^{c}_{jmn}
               =  \Ric^{g}_{ij}-\frac{\alpha^2}{2}\cdot\diag(1,\ldots,1,k).
\ee
\begin{prop}\label{prop:1}
The characteristic Ricci tensor $\Ric^c$ of an $\alpha$-Sasakian manifold vanishes
in direction of the contact vector field $\xi$, $\Ric^{c}(\xi) = 0$.
\end{prop}
\subsection{\texorpdfstring{$\alpha$}{Alpha}-Sasakian structures in dimension
\texorpdfstring{$5$}{5}.}
Fix the flux form by $F=A\cdot\ast\eta$. Using this together with
$T=B\cdot T^c$ and $\Ric^c$ leads to the first and second contractions
$K^{\nabla}(X)$, $K^{\nabla}$ relative to $\nabla=\nabla^1,\nabla^2$.
By proposition \ref{prop:1} these contractions assume a particularly
simple form along the contact vector field $\xi$,
\bea
K^{\nabla^1}(\xi) &=& -1/2\,\alpha^2\left(B^2-1\right)\cdot\eta 
                          -1/4\, A\left(B+1\right)\cdot d\eta,\\
K^{\nabla^2}(\xi) &=& -1/2\,\alpha^2\left(B^2-1\right)\cdot\eta,\\                          
K^{\nabla^1} & = & -1/2\left(\Scal^c - 3\,\alpha^2\left(B^2-1\right) + 3\,A^2\right)
                     + 1/2\,B\left(B-3\right)\cdot\left(d\eta\hut d\eta\right) \\                    
                  && + 1/2\,A\left(B-3+4\,q\right)\cdot\left(\eta\hut d\eta\right),\\
K^{\nabla^2}  &=&  -1/2\left(\Scal^c-3\,\alpha^2\left(B^2-1\right)\right)
                     + 1/2\,B\left(B-3\right)\cdot\left(d\eta\hut d\eta\right) \\
                   &&  + 2\,A\cdot\left(\eta\hut d\eta\right).                                         
\eea
As a start we consider $\nabla^1$-parallel spinor fields.
\begin{prop}
If there exists a $\nabla^1$-parallel spinor field $\Psi_1$ for $\nabla^1\neq\nabla^c$,
then the following assertions hold:
\begin{enumerate}
  \item The component of $\Psi_1$ in the one-dimensional subbundles
        defined by $\Phi\cdot\Psi=\pm2\,i\cdot\Psi$ vanishes, and the
        parameters are related by $A=\pm\alpha\,(B-1)\neq0$ respectively.
  \item The characteristic Ricci tensor is $\Ric^{c}= (B+2\,q\,(B-1))\,
        \alpha^2\cdot(g-\eta\ox\eta)$.
  \item If $\Psi_1$ contains a non-vanishing component that satisfies the
        equation $\Phi\cdot\Psi=0$, then $B=-1$ and $q=0$.
\end{enumerate}
\end{prop}\noindent
As mentioned in \autoref{sec:4} a natural question is which of these necessary
conditions are sufficient as well. We can formulate an answer for line bundles.
Let us thus assume that
\be
A=\pm\alpha\,(B-1), \quad \Ric^c = (B+2\,q\,(B-1))\,\alpha^2\cdot\diag(1,1,1,1,0).
\ee
Then the equation $K^{\nabla^1}(X)\Psi = 0$ is satisfied for all spinor
fields in the one-dimensional spin subbundle defined by $\Phi\cdot\Psi=\mp
2\,i\cdot\Psi$ (signs according), and the spinorial covariant derivative
$\nabla^1$ preserves this subbundle,
\be
\nabla_X^1\Psi-\nabla_X^c\Psi = \pm \frac{i}{2}\,\alpha\,(B-1)\,(1+2q)\,\eta(X)
\cdot\Psi\quad\mrm{for}\quad \Phi\cdot\Psi=\pm2\,i\cdot\Psi.
\ee
\begin{thm}\label{thm:1}
Let $(M^5,g,\xi,\eta,\vphi)$ be a five-dimensional, simply connected,
$\alpha$-Sasakian spin manifold with flux form $F=A\cdot\ast\eta$. If
the parameters are related by $A=\pm\alpha\,(B-1)$ and
\be
\Ric^c = (B+2\,q\,(B-1))\,\alpha^2\cdot(g-\eta\ox\eta),
\ee
then there exists a $\nabla^1$-parallel spinor field in the one-dimensional
subbundle defined by the relation $\Phi\cdot\Psi=\mp2\,i\cdot\Psi$ respectively.
\end{thm}
\begin{exa}\label{exa:1}
We provide an example in the case of $\alpha=2$, i.e.~Sasakian structures.
Simply connected, Sasakian spin manifolds with $\Ric^{c}=\diag(4,4,4,4,0)$
can be constructed -- for instance -- as bundles over four-dimensional
K\"ahler-Einstein manifolds with positive scalar curvature: Consider a simply
connected K\"ahler-Einstein manifold $(N^4,\bar{g},J)$ with scalar curvature
$\Scal^{\bar{g}}=32$. Then there exists an $S^1$-bundle $M^5\ra N^4$ with
$(M^5,g,\xi,\eta,\vphi)$ Sasakian, such that $\Ric^g=\diag(6,6,6,6,4)$ (see
\cite{FK00}). Using Theorem \ref{thm:1} we deduce the existence of a solution
$\Psi_1$ to $\nabla^1\Psi=0$ with $\Phi\cdot\Psi_1=\pm2\,i\cdot \Psi_1$ by
choosing $T=B\cdot(\eta\hut d\eta)$, $F=\mp2\,(B-1)\cdot\ast\eta$, $q=-1/2$
respectively.
\end{exa}
We proceed with the case $\nabla^2\Psi=0$.
\begin{prop}\label{prop:2}
If there exists a $\nabla^2$-parallel spinor field $\Psi_2$,
then $B=1$ and one of the following occurs:
\begin{enumerate}
  \item $\Ric^c=\alpha\,(\alpha\pm 2\,A)\cdot(g-\eta\ox\eta)$ and
        $\Phi\cdot\Psi_2=\mp 2\,i\cdot\Psi_2$.
  \item $A=0$, $\Ric^c=-\alpha^2\cdot(g-\eta\ox\eta)$ and $\Phi\cdot\Psi_2=0$.
  \item $A\neq0$, $\Scal^{c}=-4\,\alpha^2$, $\Phi\cdot\Psi_2=0$ and
        $(M^5,g)$ is not $\eta$-Einstein.
\end{enumerate}
\end{prop}
\begin{rmk}
This proposition generalizes results of \cite{FI02}.
\end{rmk}
\begin{exa}
A Sasakian structure ($\alpha=2$) that admits a Ricci tensor as in (2)
and a $\nabla^c$-parallel spinor field with $\Phi\cdot\Psi=0$, is
locally equivalent to the Sasakian structure arising from left invariant
vector fields on the five-dimensional Heisenberg group (see \cite{FI03a}).
\end{exa}\noindent
As for the converse of proposition \ref{prop:2} suppose that
\be
B=1, \quad \Ric^c=\alpha\,(\alpha\pm 2\,A)\cdot\diag(1,1,1,1,0).
\ee
Then $K^{\nabla^2}(X)\Psi=0$ is satisfied for spinor fields in the subbundle
defined by $\Phi\cdot\Psi=\mp 2\,i\cdot\Psi$ respectively. Computing the difference
\be
\nabla_X^2\Psi-\nabla_X^c\Psi = -i\, A\,\eta(X)\cdot\Psi\quad\mrm{for}
\quad \Phi\cdot\Psi=\pm2\,i\cdot\Psi,
\ee
we obtain that $\nabla^2$ preserves these spin subbundles.
\begin{thm}
Let $(M^5,g,\xi,\eta,\vphi)$ be a five-dimensional, simply connected, $\alpha$-Sasakian
spin manifold with flux form $F=A\cdot\ast\eta$. Then there exists a $\nabla^2$-parallel
spinor field in the one-dimensional subbundle defined by $\Phi\cdot\Psi=\pm2\,i\cdot\Psi$,
if $B=1$ and if the characteristic Ricci tensor is correspondingly given by
\be
\Ric^c=\alpha\,(\alpha\mp 2\,A)\cdot(g-\eta\ox\eta).
\ee
\end{thm}
\begin{rmk}
There exists no solution to $\nabla^2\Psi=0$ for $F=A\cdot\ast\eta\neq0$ on
the Sasakian structure described in example \ref{exa:1}.
\end{rmk}
\subsection{\texorpdfstring{$\alpha$}{Alpha}-Sasakian structures in dimension
\texorpdfstring{$7$}{7}.}
As in the five-dimensional case we begin with fixing the flux form,
$F=A\cdot\ast(\eta\hut\Phi)$, and then directly proceed with classification
results for the equation $\nabla^1\Psi=0$.
\begin{prop}
The existence of a $\nabla^1$-parallel spinor field $\Psi_1$ for $\nabla^1\neq\nabla^c$
leads to one of the following cases:
\begin{enumerate}
  \item The parameters are related by $A\,(4\,q-6)=\alpha\,(B-1)$ and
        $A\neq0$. The spinor field is fixed by $\Psi_1=\Psi_1^++\Psi_1^-$ and
        $\Phi\cdot\Psi_1^\pm=\pm 3\,i\cdot\Psi_1^\pm$. The characteristic
        Ricci tensor is
        \be
        \Ric^c=\alpha\,(2\,\alpha-9\,A)\cdot(g-\eta\ox\eta).
        \ee
  \item The parameters satisfy $A=-\alpha/3$ and $B=-4\,q/3-1$. The spinor field
        is given by $\Psi_1=\Psi_1^++\Psi_1^-$ and $\Phi\cdot\Psi_1^\pm=\pm
        i\cdot\Psi_1^\pm$. The characteristic Ricci tensor has the diagonal form 
        \be
        \Ric^c=\alpha^2\cdot(g-\eta\ox\eta).
        \ee
\end{enumerate}
\end{prop}\noindent
Again following \autoref{sec:4}, suppose
\be
A\,(4\,q-6)=\alpha\,(B-1), 
\quad \Ric^c=\alpha\,(2\,\alpha-9\,A)\cdot\diag(1,1,1,1,1,1,0).
\ee
Then $K^{\nabla^1}(X)\Psi = 0$ holds for all spinor fields in the
one-dimensional spin subbundles defined by $\Phi\cdot\Psi=\pm3\,i\cdot\Psi$.
The simple check of
\be
\nabla_X^1\Psi-\nabla_X^c\Psi = \mp \frac{9}{2}\,i\, A\,\eta(X)\cdot\Psi\quad
\mrm{for}\quad\Phi\cdot\Psi=\pm3\,i\cdot\Psi
\ee
allows to conclude that $\nabla^1$ preserves these subbundles.
\begin{thm}\label{thm:2}
Let $(M^7,g,\xi,\eta,\vphi)$ be a seven-dimensional, simply connected,
$\alpha$-Sasakian spin manifold with flux form $F=A\cdot\ast(\eta\hut\Phi)$.
There exist two $\nabla^1$-parallel spinor fields $\Psi_1^\pm$ satisfying
$\Phi\cdot\Psi_1^\pm=\pm 3\,i\cdot\Psi_1^\pm$, if the system parameters
are related by $A\,(4\,q-6)=\alpha\,(B-1)$ and if the characteristic Ricci
tensor is given by
\be
\Ric^c=\alpha\,(2\,\alpha-9\,A)\cdot(g-\eta\ox\eta).
\ee
\end{thm}
\begin{exa}
Simply connected Sasakian manifolds which admit the characteristic Ricci
tensor of theorem \ref{thm:2} can be constructed via the Tanno deformation
of an arbitrary seven-dimensional Einstein-Sasakian structure $(\tilde{M}^7,\tilde{g},
\tilde{\xi},\tilde{\eta},\tilde{\vphi})$. This deformation is defined by
\be
\vphi:=\tilde{\vphi},\quad \xi:=a^2\cdot\tilde{\xi},
\quad \eta:=a^{-2}\cdot\tilde{\eta},
\quad g:=a^{-2}\cdot\tilde{g}+(a^{-4}-a^{-2})\cdot\tilde{\eta}\ox\tilde{\eta}
\ee
with the deformation parameter $a^2=3/2$ (see \cite{BGN00}). We recommend
the article \cite{BGM04} for further constructions of Sasakian structures
of $\eta$-Einstein type.
\end{exa}
We close this subsection by classifying the solution space of $\nabla^2\Psi = 0$.
\begin{prop}
Suppose there exists a $\nabla^2$-parallel spinor field $\Psi_2$. Then $\Psi_2$
is parallel with respect to $\nabla^c$, $4\,A=\alpha\,(B-1)$ is satisfied and
one of the following holds:
\begin{enumerate}
	\item $\Ric^{g}=\frac{\alpha^2}{2}\cdot\diag(5,5,5,5,5,5,3)$,
	      $\Psi_0=\Psi_0^++\Psi_0^-$
	      and $\Phi\cdot\Psi_0^\pm=\pm 3\,i\cdot\Psi_0^\pm$.
  \item $\Scal^{g}=\frac{\alpha^2}{2}$, $\Psi_0=\Psi_0^++\Psi_0^-$
	      and $\Phi\cdot\Psi_0^\pm=\pm i\cdot\Psi_0^\pm$.
\end{enumerate}
In case $\mrm{(2)}$ $(M^7,g)$ is not $\eta$-Einstein.
\end{prop}
\begin{rmk}
The integrability condition (1) corresponds with \cite{FI02}*{Thm.~9.2} in
the case of seven-dimensional Sasakian structures. Due to the necessary
$\nabla^c$-parallelism of the solving spinor field we can refer to theorem
\ref{thm:2} for the construction of solutions to $\nabla^2\Psi = 0$.
\end{rmk}
%
%
%
%
%
\section{Almost Hermitian structures}\label{sec:6}\noindent
In this section we will consider almost Hermitian structures of class
$\mc{C}[\G]$ for a connected, non-abelian subgroup $\G\subset\U(3)$
that stabilizes a non-trivial $3$-form. $\G$ will be therefore up to
conjugation one of the following groups (see \cite{Sch06}):
\begin{enumerate}
  \item $\SU(3)\hookrightarrow \U(3)$.
  \item $\SO(3)\hookrightarrow \U(3),\quad$ (three-dimensional, irreducible
        complex representation of $\SU(2)$).
  \item $\SU(2)\hookrightarrow_{\iota} \U(3),\quad
        \iota:\,\,A\longmapsto\left[\begin{array}{cc}A&0\\0&1\end{array}\right]$.
  \item $\U(2)\hookrightarrow_{\iota_k} \U(3),\quad
        \iota_k:\,\,A \longmapsto \left[ \begin{array}{cc}A&0\\0&\det(A)^k
        \end{array}\right],\quad k=-1,0,1$.
\end{enumerate}
One of the main results in this section deals with the Killing spinor equation
of special type.
\begin{thm}\label{mthm:1}
Let $\G$ be a connected, non-abelian subgroup of $\U(3)$ that stabilizes a
non-trivial $3$-form and $(M^6,g,J)$ a six-dimensional, almost Hermitian spin
manifold of class $\mc{C}[\G]$ with characteristic connection $\nabla^c$,
$4$-form $F=A\cdot\ast\Omega\neq 0$, $3$-form $T=B\cdot T^c$ and spinorial
covariant derivative
\be
\nabla^{0}_{X}\Psi=\nabla^g_X\Psi + \frac{1}{4}\,(X\haken T)\cdot\Psi 
                                  + \frac{1}{2}\,(X\haken F)\cdot\Psi
                                  + (X\hut F)\cdot\Psi.
\ee
There exists a $\nabla^0$-parallel spinor field $\Psi_0$ if and only if the
following conditions are satisfied:
\begin{enumerate}
  \item The spinor field $\Psi_0$ is parallel with respect to the characteristic
        connection $\nabla^c$ and satisfies $\ast\Omega\cdot\Psi=-3\cdot\Psi$.
  \item The $3$-form $T$ coincides with the characteristic torsion, $T=T^c$.
\end{enumerate}
\end{thm}\noindent
The sufficiency of (1) and (2) follows from a direct computation after fixing
an adapted frame $(e_1,\ldots,e_6)$ and a spin representation. Due to the
complexity of the algebraic systems we have split the description into
\be
\mc{C}[\SU(3)],\quad \mc{C}[\SO(3)],\quad \mc{C}[\SU(2)]\quad \mathrm{and}\quad
\mc{C}[\U(2),{\iota_k}] 
\ee
in order to prove necessity. We shall only describe in detail the first case,
as the remaining are completely analogous.
\subsection{Almost Hermitian structures of class \texorpdfstring{$\mc{C}[\SU(3)]$}{C[SU(3)]}.}
This class is equivalent to the class of strictly (i.e.~non-K\"ahler) nearly
K\"ahler structures. The torsion form of $\nabla^c$ is given by
\be
T^c=a\cdot(-e_{246}+e_{136}+e_{145}+e_{235})\in\Lambda^3_2(\R^6),
\quad a\in\R,\,\,a>0
\ee
for a chosen adapted frame $(e_1,\ldots,e_6)$. Known formulae (see \cite{BFGK91})
\bea
\left\|(\nabla^g_XJ)(Y)\right\|^2
&=&a^2\left\{\left\|X\right\|^2\left\|Y\right\|^2-g(X,Y)^2-g(JX,Y)^2\right\},\\
(\nabla^g_XJ)(Y)&=&\frac{1}{2}\cdot\left\{J(T^c(X,Y))-T^c(X,JY)\right\}
\eea
ensure that $(M^6,g,J)$ is Einstein (see \cite{Gra76}),
\be
\Ric^g=5\,a^2\cdot g.
\ee
Then
\be
\Ric^c_{ij}=\Ric^g_{ij}-\frac{1}{4}\cdot T^c_{imn}T^c_{jmn}=4\,a^2\cdot g_{ij}.
\ee
\begin{exa}
Simply connected, homogeneous examples of strictly nearly K\"ahler structures
include the six-dimensional sphere $S^6$, the complex projective space $\C\mrm{P}
(3)$ and the flag manifold $\mrm{F}(1,2)$ (see \cite{But05}).
\end{exa}
We now study the spin geometry of the local model $\R^6$. Denote by $\Delta_6$
the space of complex spinors on $\R^6$. Vectors and forms act on $\Delta_6$ by
Clifford multiplication, for which we choose the matrix representation of \cite{Fri00}.
As an endomorphism of spinors the torsion form $T^c$ then reads 
\be
T^c=\left[\begin{array}{ccc} 0 & 4\,a\,i & 0\\
                            -4\,a\,i & 0 & 0\\
                                 0   & 0 & 0\end{array} \right].
\ee
Define the eigenspaces of $T^c$ by
\be
\Delta^{1\pm}:=\left\{\Psi\in\Delta_6\,:\,T^c\cdot\Psi=\pm4\,a\cdot\Psi\right\},
\quad\Delta^{6}:=\left\{\Psi\in\Delta_6\,:\,T^c\cdot\Psi=0\right\}. 
\ee

Consequently, the torsion form $T^c$ splits the spinor bundle $\Sigma$ of
$(M^6,g)$ into two one-dimensional subbundles $\Sigma^{1\pm}$ and one
six-dimensional subbundle $\Sigma^{6}$. We denote the components of an
arbitrary spinor field $\Psi\in\Gamma(\Sigma)$ correspondingly by
\be
\Psi=\Psi^{1+}+\Psi^{1-}+\Psi^{6}.
\ee
The ansatz for the family of $4$-forms $F$ chosen here is given by the
Hodge dual of the K\"ahler form,
\be
F:=A\cdot\ast\Omega,\quad A\in\R.
\ee

Let us consider $\nabla^1$. The spinorial covariant derivative $\nabla^1$
can be rewritten as follows:
\be
\nabla^1_X\Psi = \nabla^c_X\Psi + \frac{B-1}{4}\,(X\haken T^c)\cdot\Psi 
                                  + \frac{A}{2}\,(X\haken \ast\Omega)\cdot\Psi
                                  + A\,q\,(X\hut \ast\Omega)\cdot\Psi.
\ee
Since $T^c$ is parallel with respect to $\nabla^c$, the splitting of the
spinor bundle is preserved by $\nabla^c$. With this property we derive the
first necessary conditions on the existence of $\nabla^1$-parallel spinor
fields.
\begin{lem}\label{lem:1}
If there exists a $\nabla^1$-parallel spinor field $\Psi_1$, the following statements hold:
\begin{enumerate}
  \item If $\Psi_1\in\Gamma(\Sigma^{6})$, then $A\neq0$, $B=1$ and $q=-1$.
  \item If $\Psi_1\in\Gamma(\Sigma^{1+}\op\Sigma^{1-})$, the system parameters satisfy 
        \be
        (\ast)\quad 2\,(q-1)\,A=\pm a\,(B-1).
        \ee
        If this expression is non-trivial, then $\Psi^{1\pm}_1=0$ respectively.
  \item If $(\ast)$ holds, then $\nabla^1\Psi=\nabla^c\Psi$ is satisfied for
        arbitrary $\Psi\in\Gamma(\Sigma^{1\mp})$ respectively.
\end{enumerate}
\end{lem}
\begin{proof}
We prove (2). Fix a generic $\nabla^1$-parallel spinor field $\Psi_1\in\Gamma
(\Sigma^{1+}\op\Sigma^{1-})$ by
\be
\Psi_1=[(p_ + -p_-)\,i,(p_+ + p_-),0,0,0,0,0,0]^T.
\ee
We respectively have $\Psi_1\in\Gamma(\Sigma^{1\pm})$ if and only if the complex-valued
function $p_\mp$ vanishes. Due to the $\nabla^1$-parallelism of $\Psi_1$
\bea
-\nabla_{e_1}^{c}\Psi_1 & = & \nabla^{1}_{e_1}\Psi_1-\nabla_{e_1}^{c}\Psi_1 \\
                        & = & \frac{B-1}{4}\,(e_1\haken T^c)\cdot\Psi_1 
                               + \frac{A}{2}\,(e_1\haken \ast\Omega)\cdot\Psi_1 
                               + A\, q\,(e_1\hut \ast\Omega)\cdot\Psi_1 \\
                        & = & \left[ \begin{array}{c} 
                                     0\\
                                     0\\
                                     A\,(q-1)\,(p_+ - p_-)+\frac{a}{2}\,(B-1)\,(p_+ + p_-)\\
                                     0\\
                                     0\\
                                     \frac{a}{2}\,(B-1)\,(p_+ - p_-)\,i+A\,(q-1)\,(p_+ + p_-)\,i\\
                                     0\\
                                     0 
                                     \end{array}      \right].
\eea
Since $\nabla_{e_1}^{c}\Psi_1\in\Gamma(\Sigma^{1+}\op\Sigma^{1-})$, we deduce
\be
2\,(q-1)\,A=\pm a\,(B-1).
\ee
Writing an arbitrary vector field $X\in TM$ in components,
\be
X=X^1\cdot e_1+\ldots+X^6\cdot e_6,
\ee
we derive the following assuming $(q-1)\,A\neq 0$:
\bea
-\nabla_{X}^{c}\Psi_1 & = & \nabla^{1}_{X}\Psi_1-\nabla_{X}^{c}\Psi_1 \\
                      & = & \frac{B-1}{4}\,(X\haken T^c)\cdot\Psi_1 
                               + \frac{A}{2}\,(X\haken \ast\Omega)\cdot\Psi_1 
                               + A\, q\,(X\hut \ast\Omega)\cdot\Psi_1 \\
                      & = & a\,(B-1)\,p_{\pm}\cdot\left[ \begin{array}{c} 
                                                           0\\
                                                           0\\
                                                           (X^1 - i\, X^2)\\
                                                           (-X^3 + i\, X^4)\\
                                                           (i\, X^5 + X^6)\\
                                                           \pm(i\, X^1 - X^2)\\
                                                           \pm(-i\, X^3 + X^4)\\
                                                           \pm(X^5 + i\, X^6) 
                                                   \end{array}      \right].
\eea
We conclude that $p_\pm$ has to vanish respectively.
\end{proof}
We now start to apply the technique of \autoref{sec:4} to $\nabla^1\Psi = 0$,
and therefore assume that such a spinor field $\Psi_1$ exists. Let us consider
the three distinct cases
\begin{enumerate}
  \item $\Psi_1^{1+}+\Psi_1^{1-}\neq0$ and $\Psi_1^{6}\neq0$
  \item $\Psi_1^{1+}+\Psi_1^{1-}\neq0$ and $\Psi_1^{6}=0$
  \item $\Psi_1^{1+}+\Psi_1^{1-}=0$ and $\Psi_1^{6}\neq0$
\end{enumerate}
separately. Recall $\Psi_1$ is an element in the kernel of $K^{\nabla^1}(e_i)$
and $K^{\nabla^1}$,
\bea
(\ast) && K^{\nabla^1}(e_i)\Psi_1=0,\quad K^{\nabla^1}\Psi_1=0, \quad i=1,\ldots,6.
\eea
\subsubsection*{Case (1)}
We use the endomorphisms $K^{\nabla^1}$, $K^{\nabla^1}(e_6)$, $K^{\nabla^1}(e_4)$
and $K^{\nabla^1}(e_2)$. In the chosen representation these are given by
\be
K^{\nabla^1}=\left[
\begin{array}{cccccccc}
	-6\, m_1 & 6\,i\, m_2 & 0 & 0 & 0 & 0 & 0 & 0 \\
	-6\,i\, m_2 & -6\, m_1 & 0 & 0 & 0 & 0 & 0 & 0 \\
	0 & 0 & c_1 & 0 & 0 & 0 & 0 & 0 \\
	0 & 0 & 0 & c_1 & 0 & 0 & 0 & 0 \\
	0 & 0 & 0 & 0 & c_1 & 0 & 0 & 0 \\
	0 & 0 & 0 & 0 & 0 & c_1 & 0 & 0 \\
	0 & 0 & 0 & 0 & 0 & 0 & c_1 & 0 \\
	0 & 0 & 0 & 0 & 0 & 0 & 0 & c_1 
\end{array} \right],
\ee
\be
K^{\nabla^1}(e_6)=\left[
\begin{array}{cccccccc}
	0 & 0 & 0 & 0 & i\, n_1 & 0 & 0 & n_2 \\
	0 & 0 & 0 & 0 & n_2 & 0 & 0 & -i\, n_1 \\
	0 & 0 & 0 & 0 & 0 & 0 & -i\, c_2 & 0 \\
	0 & 0 & 0 & 0 & 0 & i\, c_2 & 0 & 0 \\
	i \, m_1 & m_2 & 0 & 0 & 0 & 0 & 0 & 0 \\
	0 & 0 & 0 & i\, c_2 & 0 & 0 & 0 & 0 \\
	0 & 0 & -i\, c_2 & 0 & 0 & 0 & 0 & 0 \\
	m_2 & -i\, m_1 & 0 & 0 & 0 & 0 & 0 & 0 
\end{array} \right],
\ee
\be
K^{\nabla^1}(e_4)=\left[
\begin{array}{cccccccc}
	0 & 0 & 0 & n_1 & 0 & 0 & i \, n_2 & 0 \\
	0 & 0 & 0 & -i\, n_2 & 0 & 0 & n_1 & 0 \\
	0 & 0 & 0 & 0 & 0 & 0 & 0 & c_2 \\
	-m_1 & i\, m_2 & 0 & 0 & 0 & 0 & 0 & 0 \\
	0 & 0 & 0 & 0 & 0 & -c_2 & 0 & 0 \\
	0 & 0 & 0 & 0 & c_2 & 0 & 0 & 0 \\
	-i\, m_2 & -m_1 & 0 & 0 & 0 & 0 & 0 & 0 \\
	0 & 0 & -c_2 & 0 & 0 & 0 & 0 & 0 
\end{array} \right],
\ee
\be
K^{\nabla^1}(e_2)=\left[
\begin{array}{cccccccc}
	0 & 0 & -n_1 & 0 & -i\, n_2 & 0 & 0 & 0 \\
	0 & 0 & i\, n_2 & 0 & -n_1 & 0 & 0 & 0 \\
	m_1 & -i\, m_2 & 0 & 0 & 0 & 0 & 0 & 0 \\
	0 & 0 & 0 & 0 & 0 & 0 & 0 & c_2 \\
	0 & 0 & 0 & 0 & 0 & 0 & -c_2 & 0 \\
	i\, m_2 & m_1 & 0 & 0 & 0 & 0 & 0 & 0 \\
	0 & 0 & 0 & 0 & c_2 & 0 & 0 & 0 \\
	0 & 0 & 0 & -c_2 & 0 & 0 & 0 & 0 
\end{array} \right].
\ee
Here the numbers $m_1$, $m_2$, $n_1$, $n_2$, $c_1$ and $c_2$ depend on the
system parameters in the following way:
\be
m_1= a^2/2\left(B-1\right)\left(B-5\right)+2\,A^2\left(q-1\right)\left(3\,q-1\right),
\,\, m_2=2\,a\,A\left(B-3-2\,q\left(B-2\right)\right),
\ee
\be
n_1= a^2/2\left(5-3\,B\left(B-2\right)\right)+2\,A^2\left(3-5\,q^2\right),
\,\, n_2=-4\,q\,a\,A\,(B-2),
\ee
\be
c_1= a^2\left(B\left(5\,B-6\right)-15\right)+4\,A^2\left(q\left(7\,q-4\right)-5\right),
\ee
\be
c_2= a^2/2\left(B^5-5\right)+2\,A^2\left(q\left(q-2\right)-1\right).
\ee
Forcing the corresponding matrices to have non-maximal rank yields
\be
m_1^2=m_2^2,\quad n_1^2=n_2^2,\quad c_1=0,\quad c_2=0,\\
\ee
in contradiction to the assumption.
\subsubsection*{Case (2)}
Lemma \ref{lem:1} leads to one of the following cases:
\begin{itemize}
 \item $A=0$, $B=1$ and $\Psi_1\in\Gamma(\Sigma^{1+}\op\Sigma^{1-})$.
 \item $A\neq0$, $B=1$, $q=1$ and $\Psi_1\in\Gamma(\Sigma^{1+}\op\Sigma^{1-})$.
 \item $2\,A\,(q-1)=a\,(B-1)\neq0$ and $\Psi_1\in\Gamma(\Sigma^{1-})$.
 \item $2\,A\,(q-1)=-a\,(B-1)\neq0$ and $\Psi_1\in\Gamma(\Sigma^{1+})$.
\end{itemize}
$\nabla^c\Psi_1=0$ holds and system $(\ast)$ is satisfied for all of them.
\subsubsection*{Case (3)}
Lemma \ref{lem:1} reduces the problem to one of the two cases
\begin{itemize}
  \item $A=0$, $B=1$.
  \item $A\neq0$, $B=1$, $q=-1$.
\end{itemize}
In analogy to case (1) the system $n_1^2=n_2^2$, $c_1=0$, $c_2=0$
becomes inconsistent.

Summarizing the cases (1)--(3) and using lemma \ref{lem:1} proves the
following theorem.
\begin{thm}
Let $(M^6,g,J)$ be a six-dimensional, almost Hermitian spin manifold of
class $\mc{C}[\SU(3)]$ with flux form $F=A\cdot\ast\Omega$. Then there
exists a $\nabla^1$-parallel spinor field $\Psi_1$ if and only if the
following conditions are satisfied:
\begin{enumerate}
  \item The spinor field $\Psi_1$ is parallel with respect to $\nabla^c$
        and satisfies $\ast\Omega\cdot\Psi=-3\cdot\Psi$. 
  \item The system parameters satisfy $2\,(q-1)\,A=\pm a\,(B-1)$. If this
        expression is non-trivial, then $T^c\cdot\Psi_1=\mp4\,a\cdot\Psi_1$
        respectively.
\end{enumerate}
\end{thm}
\begin{rmk}
If we fix the spinorial field equation $\nabla^1\Psi=0$ by requiring $q=1$, i.e.~if we
consider the special type $\nabla^0\Psi=0$, the $3$-form $T$ coincides with $T^c$ and
we have proven theorem \ref{mthm:1} for structures of class $\mc{C}[\SU(3)]$.
\end{rmk}
Solving $\nabla^1\Psi=0$ reduces to $\nabla^c\Psi=0$. A direct computation
proves the following.
\begin{prop}\label{prop:4}
On almost Hermitian spin manifolds of class $\mc{C}[\SU(3)]$ there exist
two $\nabla^c$-parallel spinor fields $\Psi_\pm$ such that $T^c\cdot\Psi_\pm
=\pm4\,a\cdot\Psi_\pm$.
\end{prop}
We conclude with $\nabla^2\Psi=0$.
\begin{prop}
If there exists a $\nabla^2$-parallel spinor field $\Psi_2$ with
$\nabla^2\neq\nabla^c$, the following conditions are satisfied:
\begin{enumerate}
  \item The system parameters fulfill $2\,A=\pm a\,(B-1)$ and the
        component of $\Psi_2$ in the one-dimensional spin subbundle
        defined by $T^c\cdot\Psi=\pm4\,a\cdot\Psi$ respectively, vanishes.
  \item If $B\neq2$, the spinor field $\Psi_2$ is $\nabla^c$-parallel
        and fixed by $\ast\Omega\cdot\Psi=-3\cdot\Psi$.
\end{enumerate}
\end{prop}
\subsection{Almost Hermitian structures of class \texorpdfstring{$\mc{C}[\SO(3)]$}{C[SO(3)]}.}
The characteristic torsion splits into two components,
\be
T^c = T^c_2 + T^c_{12}\in\Lambda^3_2(\R^6)\op\Lambda^3_{12}(\R^6).
\ee
We fix an adapted frame such that
\bea
T^c_2 &=& a\cdot(-e_{135}+e_{146}+e_{236}+e_{245})\in\Lambda^3_2(\R^6),\\
T^c_{12} &=& (b+c\, J)\cdot(3\, e_{135}+e_{146}+e_{236}+e_{245})\in\Lambda^3_{12}(\R^6)
\eea
for real parameters $a,b,c$. In this the characteristic Ricci tensor
is
\be
\Ric^c=4\,(a^2-b^2-c^2)\cdot g.
\ee
\begin{exa}
Homogeneous examples ($M^6 = \G/\mrm{H}$) of class $\mc{C}[\SO(3)]$ can be constructed
using one of the following spaces as the base manifold $M^6$: $\SL(2,\C)$, $S^3\x S^3$,
$\tilde{E}_3$, $N$ (see \cite{Sch06}). Here $\tilde{E}_3=\SU(2)\ltimes\R^3$ is the
universal covering of the group of Euclidean motions of $\R^3$ and $N$ a nilpotent Lie
group. 
\end{exa}
We then come to results regarding $\nabla^0\Psi=0$.
\begin{prop}
Assuming $F=A\cdot\ast\Omega\neq0$ there exists a $\nabla^0$-parallel
spinor field $\Psi_0$ if and only if the following conditions are satisfied:
\begin{enumerate}
  \item The spinor field $\Psi_0$ is parallel with respect to $\nabla^c$ and
        satisfies $\ast\Omega\cdot\Psi=-3\cdot\Psi$.
  \item The parameter $B$ equals $1$.
\end{enumerate}
If $A=0$ instead, $\nabla^0\Psi_0=0$ implies either
\begin{enumerate}
  \item $\Ric^c\neq0$, $B=1$ and $\ast\Omega\cdot\Psi_0=-3\cdot\Psi_0$ or 
  \item $\Ric^c=0$ and $B=\pm1$.
\end{enumerate}
\end{prop}
A direct computation leads to the existence of two solutions to $\nabla^c\Psi=0$
which are eigenspinors of $T^c$, 
\be
T^c\cdot\Psi_\pm = \pm 2\,\left\|T^c_2\right\|\cdot \Psi_\pm.
\ee
\begin{prop}
Almost Hermitian spin manifolds of class $\mc{C}[\SO(3)]$ admit two
$\nabla^c$-parallel spinor fields $\Psi_\pm$ satisfying $\ast\Omega
\cdot\Psi=-3\cdot\Psi$.
\end{prop}
We conclude by constructing solutions to the non-special type of Killing
spinor equation based on these $\nabla^c$-parallel spinor fields.
\begin{thm}
Let $(M^6,g,J)$ be a six-dimensional, almost Hermitian spin manifold of
class $\mc{C}[\SO(3)]$ with flux form $F=A\cdot\ast\Omega$. The equation
\bea
(\ast)&& s\,(X\haken T^c)\cdot\Psi_0 + p\,(X\haken F)\cdot\Psi_0 + q(X\hut F)\cdot\Psi_0=0
\eea
holds for the ansatz $\Psi_0=p_+\cdot\Psi_+ +p_-\cdot\Psi_-$ if and only if
it is solved by $p_+\cdot\Psi_+$ and $p_-\cdot\Psi_-$ separately, and if
$b\cdot s=c\cdot s=0$. The spinor field $\Psi_\pm$ is a solution to $(\ast)$
if and only if
\be
(2\,p-q)\,A=\pm2\,a\,s.
\ee
\end{thm}
\begin{rmk}
$\nabla^c$-parallel spinor fields $\Psi_\pm$ solving the Killing spinor equation
for $2p\neq q$ (non-special type) and $A\neq0$ force the almost Hermitian structure
$(M^6,g,J)$ to belong to the class $\mc{C}[\SU(3)]$, i.e.~$b=c=0$.
\end{rmk}
\subsection{Almost Hermitian structures of class \texorpdfstring{$\mc{C}[\SU(2)]$}{C[SU(2)]}.}
The torsion form $T^c$ splits into two components,
\be
T^c=T^c_{14}+T^c_{6}\in\Lambda^3_2(\R^6)\op\Lambda^3_{12}(\R^6)\op\Lambda^3_6(\R^6),
\ee
whose expression in an adapted frame reads
\bea
T^c_{14}&=&a\cdot(e_{145}+e_{235})\in\Lambda^3_{2}(\R^6)\op\Lambda^3_{12}(\R^6),\\
T^c_{6}&=&b\cdot(e_{125}+e_{345})\in\Lambda^3_6(\R^6)
\eea
for some real parameters $a,b$. In this frame we define the $2$-forms
$\Omega_1$ and $\Omega_2$ which are $\nabla^c$-parallel for structures
of class $\mc{C}[\SU(2)]$ and $\mc{C}[\U(2),\iota_k]$,
\be
\Omega_1:=e_{56},\quad \Omega_2:=e_{12}+e_{34}.
\ee
The characteristic Ricci tensor $\Ric^c$ is given by
\be
\Ric^c=(a^2+b^2)\cdot \diag(1,1,1,1,0,0).
\ee
\begin{exa}
The Hopf fibration $S^1\ra S^5\ra \C\mrm{P}(2)$ gives rise to a Sasakian structure
on the five-dimensional sphere $S^5$. There exists an $S^1$-bundle $M^6\ra S^5$ carrying
an almost Hermitian structure of type $\W_1\op\W_3\op\W_4$ such that the corresponding
torsion form $T^c$ is parallel with respect to $\nabla^c$ and
\be
\Hol(\nabla^c)\subset\SU(2)\subset\Iso(T^c)
\ee
is satisfied (see \cite{Sch07}). Every torsion form, i.e.~every combination of the
parameters $a$ and $b$ can be realized using this method of construction.
\end{exa}
We directly come to results regarding the Killing spinor equation of special type.
\begin{thm}
Let $(M^6,g,J)$ be a six-dimensional, almost Hermitian spin manifold of class
$\mc{C}[\SU(2)]$ with flux form $F=A_1\cdot\ast\Omega_1 + A_2\cdot\ast\Omega_2$.
Assuming $B\neq 2$ there exists a $\nabla^0$-parallel spinor field $\Psi_0$ if
and only if the following conditions are satisfied:
\begin{enumerate}
  \item The spinor field $\Psi_0$ is parallel with respect to $\nabla^c$ and
        satisfies $\ast\Omega_1\cdot\Psi=-\Psi$.
  \item The parameter $B$ equals $1$.
  \item One of the following occurs:
        \begin{itemize}
          \item $A_1=A_2\neq0$ and $\ast\Omega_2\cdot\Psi_0=-2\cdot\Psi_0$.
          \item $A_1=-A_2\neq0$ and $\ast\Omega_2\cdot\Psi_0=2\cdot\Psi_0$.
          \item $A_1=A_2=0$.
        \end{itemize}
\end{enumerate}
If $B=2$ instead, $\nabla^0\Psi_0=0$ implies $A_1^2\neq A_2^2$ and
$\ast\Omega_1\cdot\Psi_0=-\Psi_0$.
\end{thm}
A direct computation yields the following.
\begin{prop}\label{prop:3}
Any almost Hermitian structure of class $\mc{C}[\SU(2)]$ admits four
$\nabla^c$-parallel spinor fields $\Psi_1^{\pm}$, $\Psi_2^{\pm}$ such
that $\ast\Omega_1\cdot\Psi_i^{\pm}=-\Psi_i^{\pm}$ and $\ast\Omega_2
\cdot\Psi_i^{\pm}=\pm2\cdot\Psi_i^{\pm}$.
\end{prop}
We then answer the question which Killing spinor equations of non-special type
can be solved by these spinor fields.
\begin{thm}\label{thm:3}
Let $(M^6,g,J)$ be a six-dimensional, almost Hermitian spin manifold of class
$\mc{C}[\SU(2)]$ with flux form $F=A_1\cdot\ast\Omega_1 + A_2\cdot\ast\Omega_2$.
The equation
\bea
s\,(X\haken T^c)\cdot\Psi_0 + p\,(X\haken F)\cdot\Psi_0 + q(X\hut F)\cdot\Psi_0=0
\eea
is solved by
\be
\Psi_0=p_1^+\cdot\Psi_1^++p_2^+\cdot\Psi_2^++p_1^-\cdot\Psi_1^-+p_2^-\cdot\Psi_2^-
\ee
if and only if the torsion term vanishes ($s=0$) and if one of the following
occurs:
\begin{enumerate}
  \item $p_i^+=0$, $A_1=A_2$ and $2\,p=q$.
  \item $p_i^+=0$, $A_1=-2\,A_2$ and $p=-q$.
  \item $p_i^-=0$, $A_1=-A_2$ and $2\,p=q$.
  \item $p_i^-=0$, $A_1=2\,A_2$ and $p=-q$. 
\end{enumerate}
\end{thm}
\subsection{Almost Hermitian structures of class \texorpdfstring{$\mc{C}[\U(2),\iota_0]$}{C[U(2),0]}.}
There exists an adapted frame such that the characteristic torsion is given by
\be
T^c=a\cdot(e_{125}+e_{345})\in\Lambda^3_6(\R^6)
\ee
for a positive real parameter $a$. The characteristic Ricci tensor is
\be
\Ric^c=\left[
\begin{array}{cccccc}
	U_1+U_2 & 0 & V_1 & V_2 & 0 & 0 \\
	0 & U_1+U_2 & -V_2 & V_1 & 0 & 0\\
	V_1 & -V_2 & U_2 & 0 & 0 & 0 \\
	V_2 & V_1 & 0 & U_2 & 0 & 0 \\
	0 & 0 & 0 & 0 & 0 & 0 \\
	0 & 0 & 0 & 0 & 0 & 0 
\end{array} \right]
\ee
where $U_i,\,V_i\in\mc{C}^\infty(M)$ are smooth functions.
\begin{thm}
Let $(M^6,g,J)$ be a six-dimensional, almost Hermitian spin manifold of class
$\mc{C}[\U(2),\iota_0]$ with flux form $F=A_1\cdot\ast\Omega_1 + A_2\cdot\ast
\Omega_2$. Assuming $\nabla^0\neq\nabla^c$ there exists a $\nabla^0$-parallel
spinor field $\Psi_0$ if and only if the following conditions are satisfied:
\begin{enumerate}
  \item The spinor field $\Psi_0$ is parallel with respect to $\nabla^c$.
  \item The $3$-form $T$ coincides with $T^c$.
  \item One of the following occurs:
        \begin{itemize}
          \item $A_1=+A_2\neq0$ and $\ast\Omega_2\cdot\Psi_0=-2\cdot\Psi_0$.
          \item $A_1=-A_2\neq0$ and $\ast\Omega_2\cdot\Psi_0=2\cdot\Psi_0$.
        \end{itemize}
\end{enumerate}
\end{thm}
A direct computation now yields the following.
\begin{prop}
Almost Hermitian structures of class $\mc{C}[\U(2),\iota_0]$ that admit
a $\nabla^c$-parallel spinor field belong to the class $\mc{C}[\SU(2)]$. 
\end{prop}
\begin{rmk}
The existence of spinor fields solving $\nabla^0\Psi=0$ or $\nabla^c\Psi=0$
forces the almost Hermitian structure to belong to the class $\mc{C}[\SU(2)]$.
We can therefore refer to proposition \ref{prop:3} and theorem \ref{thm:3} for
the construction of solutions to the Killing spinor equation of non-special type.
\end{rmk}
\subsection{Almost Hermitian structures of class \texorpdfstring{$\mc{C}[\U(2),\iota_1]$}{C[U(2),1]}.}
Using an adapted frame the torsion form can be written as
\be
T^c=a\cdot(e_{135}-e_{245}+e_{236}+e_{146})\in\Lambda^3_{12}(\R^6).
\ee
for some positive real parameter $a$. The characteristic Ricci tensor
is proportional to the metric, $\Ric^c=4\,a^2\cdot g$.

We consider the Killing spinor equation of special type.
\begin{thm}
Let $(M^6,g,J)$ be a six-dimensional, almost Hermitian spin manifold of class
$\mc{C}[\U(2),\iota_1]$ with flux form $F=A_1\cdot\ast\Omega_1 + A_2\cdot\ast\Omega_2$.
Then there exists a $\nabla^0$-parallel spinor field $\Psi_0$ if and only if the
following conditions are satisfied:
\begin{enumerate}
  \item The spinor field $\Psi_0$ is parallel with respect to $\nabla^c$ and
        fixed by $\ast\Omega_2\cdot\Psi_0=2\cdot\Psi_0$.
  \item The parameter $B$ equals $1$.
  \item The flux form parameters satisfy $A_1=-A_2$.    
\end{enumerate}
\end{thm}
A direct computation yields the following on the existence of $\nabla^c$-parallel
spinor fields.
\begin{prop}
There exist no $\nabla^c$-parallel spinor fields for almost Hermitian structures
of class $\mc{C}[\U(2),\iota_1]$.
\end{prop}
\subsection{Almost Hermitian structures of class \texorpdfstring{$\mc{C}[\U(2),\iota_{-1}]$}{C[U(2),-1]}.}
This class is a subclass of $\mc{C}[\SU(3)]$.
\begin{exa}
An almost Hermitian structure of this class is locally isomorphic to the
complex projective space $\C\mrm{P}(3)$ equipped with the nearly K\"ahler
structure coming from the twistor construction, realized by $\SO(5)/\U(2)$
(see \cite{Sch06}).
\end{exa}
We then directly state results regarding $\nabla^0\Psi=0$ using a different
ansatz on $F$ than for the class $\mc{C}[\SU(3)]$.
\begin{thm}
Let $(M^6,g,J)$ be a six-dimensional, almost Hermitian spin manifold of class
$\mc{C}[\U(2),\iota_{-1}]$ with flux form $F=A_1\cdot\ast\Omega_1 + A_2\cdot
\ast\Omega_2$. Then there exists a $\nabla^0$-parallel spinor field $\Psi_0$
if and only if the following conditions are satisfied:
\begin{enumerate}
  \item The spinor field $\Psi_0$ is parallel with respect to $\nabla^c$ and
        fixed by $\ast\Omega_2\cdot\Psi_0=-2\cdot\Psi_0$.
  \item The parameter $B$ takes the value $1$.
  \item The flux form parameters satisfy $A_1=A_2$.    
\end{enumerate}
\end{thm}
\begin{rmk}
The existence of $\nabla^0$-parallel spinor fields forces the $4$-form to be
proportional to $\ast\Omega$.
\end{rmk}
We conclude by constructing solutions to the Killing spinor equation of non-special
type using the $\nabla^c$-parallel spinor fields $\Psi_\pm$ of proposition \ref{prop:4}.
\begin{thm}
Let $(M^6,g,J)$ be a six-dimensional, almost Hermitian spin manifold of class
$\mc{C}[\U(2),\iota_{-1}]$ with flux form $F=A_1\cdot\ast\Omega_1 + A_2\cdot
\ast\Omega_2$. The equation
\bea
(\ast)&&s\,(X\haken T^c)\cdot\Psi_0 + p\,(X\haken F)\cdot\Psi_0 + q(X\hut F)\cdot\Psi_0=0
\eea
holds for ansatz $\Psi_0=p_+\cdot\Psi_++p_-\cdot\Psi_-$ if and only if it
is solved by $p_+\cdot\Psi_+$ and $p_-\cdot\Psi_-$ separately. The spinor field
$\Psi_\pm$ is a solution to $(\ast)$ if and only if
\be
(p+q)\,(A_1-A_2)=0,\quad q\,A_1-2\,p\,A_2\pm2\,a\,s=0.
\ee
\end{thm}
\begin{rmk}
Qualitatively, there exist two possible deformations of the equation
$\nabla^c\Psi=0$ leading to a Killing spinor equation with non-vanishing
flux form such that either $\Psi_+$ or $\Psi_-$ is a solution.
\end{rmk}
%
%
%
%
\section{Cocalibrated \texorpdfstring{$\G_2$}{G2}-structures}\label{sec:7}\noindent
In this section we study nearly parallel $\G_2$-structures and cocalibrated
$\G_2$-structures of class $\mc{C}[\g]$ where $\g$ is a proper, non-abelian
subalgebra of $\g_2$. There exist up to conjugation eight subalgebrae of this
type \cite{Dyn57},
\bea
\,&\su(3),\quad\so(3)\subset\su(3),\quad\su(2)\subset\su(3),\quad\un(2)\subset\su(3),&\\
\,&\su_c(2),\quad\R\op\su_c(2),\quad\su(2)\op\su_c(2),\quad\so_{ir}(3).&
\eea
Cocalibrated $\G_2$-structures of class $\mc{C}[\so_{ir}(3)]$ do not exist
(see \cite{Fri06}). Let us state the main result concerning the Killing
spinor equation of special type.
\begin{thm}\label{mthm:2}
Let $\g$ be a proper, non-abelian subalgebra of $\g_2$ and $(M^7,g,\omega^3)$ a
cocalibrated $\G_2$-manifold of class $\mc{C}[\g]$ with characteristic connection
$\nabla^c$, characteristic torsion $T^c$, $4$-form $F=A\cdot\ast\omega^3\neq 0$,
$3$-form $T=B\cdot T^c$ and spinorial covariant derivative
\be
\nabla^{0}_{X}\Psi=\nabla^g_X\Psi + \frac{1}{4}\,(X\haken T)\cdot\Psi 
                                  + \frac{3}{4}\,(X\haken F)\cdot\Psi
                                  + (X\hut F)\cdot\Psi.
\ee
In case $B\neq-7$, there exists a $\nabla^0$-parallel spinor field $\Psi_0$
if and only if the following conditions are satisfied:
\begin{enumerate}
  \item The spinor field $\Psi_0$ is parallel with respect to $\nabla^c$ and
        satisfies $\ast\omega^3\cdot\Psi=-7\cdot\Psi$.
  \item The $3$-form $T$ coincides with the torsion form, $T=T^c$.
\end{enumerate}
\end{thm}\noindent
The sufficiency of (1) and (2) can be checked in a direct computation after
fixing an adapted frame $(e_1,\ldots,e_7)$ and a spin representation. For the
necessity we split the consideration into $\mc{C}[\su(3)]\cup\mc{C}[\so(3)]$,
$\mc{C}[\su(2)]\cup\mc{C}[\un(2)]$ and $\mc{C}[\su_c(2)]\cup\mc{C}[\R\op\su_c(2)]
\cup\mc{C}[\su(2)\op\su_c(2)]$. The latter will be called $\mc{C}[\su_c(2)\,\mrm{rel.}]$
henceforth.
\subsection{Nearly parallel \texorpdfstring{$\G_2$}{G2}-structures.}
Recall that the fundamental form $\omega^3$ of a given nearly parallel
$\G_2$-structure $(M^7,g,\omega^3)$ satisfies $d\omega^3=-\lambda\cdot
\ast\omega^3$ for a real parameter $\lambda\neq0$ (see \autoref{sec:3}).
The torsion form is given by
\be
T^c=-\frac{\lambda}{6}\cdot\omega^3\in\Lambda^3_1(\R^7).
\ee
Obviously, $T^c$ is parallel with respect to $\nabla^c$. The characteristic
Ricci tensor has the diagonal form
\be
\Ric^c=\frac{\lambda^2}{3}\cdot g.
\ee

Let us study the spin geometry of nearly parallel $\G_2$-structures. The existence
of a $3$-form satisfying the differential equation $d\omega^3=-\lambda\cdot\ast
\omega^3$ is equivalent to the existence of a real Killing spinor (see \cite{FKMS97}),
$\nabla^g_X\Psi=\frac{\lambda}{8}\, X\cdot\Psi$ . Denote by $KS(M^7,g)$ the space
of all Killing spinors to the Killing number $\lambda/8$,
\be
KS(M^7,g):=\{\Psi\in\Gamma(\Sigma) \, : \, \nabla^g_X\Psi=(\lambda/8)\, X\cdot\Psi\}.
\ee
The dimension of $KS(M^7,g)$ is bounded by three in the case of compact, simply connected
$M^7\neq S^7$ (see \cite{BFGK91}), $1\leq\dim[KS(M^7,g)]\leq3$. We say that a nearly parallel
$\G_2$-structure \textit{is of type} $i$, if the dimension of $KS(M^7,g)$ equals $i$.
\begin{exa}
Nearly parallel $\G_2$-structures of type $3$ are $3$-Sasakian manifolds (see \cite{FK90}).
The only regular examples are $S^7$ and $N(1,1)=\SU(3)/S^1$ (see \cite{BFGK91}). Those of
type $2$ are Einstein-Sasakian manifolds (see \cite{FK90}) and can be constructed as
circle bundles over six-dimensional K\"ahler-Einstein manifolds $X^6$. A homogeneous example
is $N(1,1)$ with $X^6=F(1,2)$. Finally, an example of type $1$ is $\SO(5)/\SO(3)$ equipped
with a special Riemannian metric of Bryant (see \cite{Bry87}).
\end{exa}\noindent
As an endomorphism of spinors the fundamental form $\omega^3$ reads
\be
\omega^3=\diag(-7,1,1,1,1,1,1,1).
\ee
We split an arbitrary spinor field $\Psi$ into two components,
\be
\Psi=\Psi^1+\Psi^7,\quad \omega^3\cdot\Psi^1=-7\cdot\Psi^1,
\quad\omega^3\cdot\Psi^7=\Psi^7.
\ee
\begin{prop}\label{prop:5}
Let $\Psi_0\in KS(M^7,g)$ be a real Killing spinor. Then $\Psi_0^1,\Psi_0^7\in KS(M^7,g)$.
\end{prop}
\begin{proof}
A direct computation yields
\be
\frac{\lambda}{8}\, X\cdot\bar{\Psi}+\frac{1}{4}\,(X\haken T^c)\cdot\bar{\Psi}=0,\quad
\left(\frac{\lambda}{8}\, X\cdot\hat{\Psi}+\frac{1}{4}\,(X\haken T^c)\cdot\hat{\Psi}
\right)^1=0
\ee
for arbitrary spinor fields $\bar{\Psi}$ and $\hat{\Psi}$ satisfying $\omega^3\cdot\Psi
=-7\cdot\Psi$ and $\omega^3\cdot\Psi=\Psi$, respectively. A real Killing spinor $\Psi_0
\in KS(M^7,g)$ therefore satisfies $(\nabla^c_X\Psi)^1=0$ and hence $\Psi_0^1$ is
$\nabla^c$-parallel.
\end{proof}
We now move on to classification results relative to $\nabla^1\Psi = 0$.
\begin{prop}\label{prop:6}
Let $F=A\cdot\ast\omega^3$ be the flux form. If there exists a $\nabla^1$-parallel
spinor field $\Psi_1$, one of the following holds:
\begin{enumerate}
  \item The component $\Psi_1^7$ vanishes. $\Psi_1$ is a real Killing spinor,
        $\Psi_1\in KS(M^7,g)$. The parameters satisfy $-24\,A\,(q-1)=
        \lambda\,(B-1)$.
  \item The component $\Psi_1^1$ vanishes. $\Psi_1$ is a real Killing spinor,
        $\Psi_1\in KS(M^7,g)$. The parameters satisfy $A=\lambda/6$,
        $B=-4\,q-3$.
  \item Both components of $\Psi_1$ are non-trivial, and the parameters satisfy either
        \begin{itemize}
        \item $A=\lambda/3$ and $B=-8\,q+9$ or
        \item $A=0$ and $B=3$.
        \end{itemize}
\end{enumerate}
\end{prop}
A short computation directly leads to the following.
\begin{thm}\label{thm:4}
Let $(M^7,g,\omega^3)$ be a nearly parallel $\G_2$-manifold with flux form
$F=A\cdot\ast\omega^3$. A spinor field $\Psi_1$ is $\nabla^1$-parallel,
if it is a real Killing spinor, $\Psi_1\in KS(M^7,g)$, and if one of the
following holds:
\begin{itemize}
\item $-24\,A\,(q-1)=\lambda\,(B-1)$ and $\Psi_1^7=0$.
\item $A=\lambda/6$, $B=-4\,q-3$ and $\Psi_1^1=0$.
\end{itemize}
\end{thm}
\begin{rmk}
There exist $i$ Killing spinors for nearly parallel $\G_2$-structures of type
$i$. Without loss of generality, each of these satisfies either $\Psi_1^7=0$ or
$\Psi_1^1=0$ (cf.~proposition \ref{prop:5}).
\end{rmk}
\begin{rmk}
Let us compare the previous results to those of \cite{AF03}. We therefore consider
a simply connected, nearly parallel $\G_2$-structure, and normalize its metric such
that the scalar curvature equals $168$. Given this there exists a real Killing spinor
$\Psi$ ($\nabla^g_X\Psi = X\cdot\Psi$) with $\omega^3\cdot\Psi=-7\cdot\Psi$. The
equation
\be
X\cdot\Psi+\frac{r}{4}\,(X\haken\omega^3)\cdot\Psi+s\,(X\haken\ast\omega^3)\cdot\Psi
+t\,(X\hut\ast\omega^3)\cdot\Psi = 0
\ee
holds for this spinor field if and only if $16\,s=-4+12\,q-3\,r$ (see \cite{AF03}).
If we translate the latter relation into the conventions of this paper,
\be
r = -\frac{4}{3}\, B, \quad s = \frac{3}{4}\, A, \quad t = q\, A,
\ee
we obtain $-24\,A\,(q-1)=8\,(B-1)$, and the result above corresponds with
proposition \ref{prop:6} and theorem \ref{thm:4}.
\end{rmk}
We proceed with classification results regarding $\nabla^2\Psi=0$.
\begin{prop}
Let $F=A\cdot\ast\omega^3$ be the flux form. If there exists a $\nabla^2$-parallel
spinor field $\Psi_2$, one of the following assertions hold:
\begin{enumerate}
  \item The parameters satisfy $-24\,A=\lambda\,(B-1)$. $\Psi_2$ is a real
        Killing-Spinor, $\Psi_2\in KS(M^7,g)$, and its component $\Psi_2^7$ vanishes.
  \item The parameters satisfy $-24\,A=\lambda\,(B-3)$. Both components of
        $\Psi_2$ are non-trivial, i.e.~$\Psi_2^1,\Psi_2^7\neq0$. 
\end{enumerate}
\end{prop}
We conclude with results on the existence of $\nabla^2$-parallel spinor fields.
\begin{thm}
Let $(M^7,g,\omega^3)$ a nearly parallel $\G_2$-manifold with flux form
$F=A\cdot\ast\omega^3$. A spinor field $\Psi_2$ is parallel with respect
to $\nabla^2$, if it is a real Killing spinor, $\Psi_2\in KS(M^7,g)$, with
$\Psi_2^7=0$, and if the relation $-24\,A=\lambda\,(B-1)$ holds.
\end{thm}
\subsection{Cocalibrated \texorpdfstring{$\G_2$}{G2}-structures of class
\texorpdfstring{$\mc{C}[\su(3)]\cup\mc{C}[\so(3)]$}{C[su(3)], C[so(3)]}.}\label{ssec:1}
There are two different types of admissible torsion forms. In an adapted frame these
are given by
\bea
T^c_\mrm{I} &=& a\cdot( e_{127}+e_{347}+e_{567} ),\\
T^c_\mrm{II}&=& a\cdot(-2\,e_{123}+e_{136}-e_{145}+e_{235}+e_{246}+2\,e_{356})\\
         &&+b\cdot(-2\,e_{124}-e_{135}-e_{146}+e_{236}-e_{245}+2\,e_{456})\\
         &&+c\cdot( e_{135}-e_{146}-e_{236}-e_{245} )
\eea
where $a,b,c$ are real parameters. The characteristic Ricci tensor has the diagonal
form
\be
\Ric^c=\lambda\cdot \diag(1,1,1,1,1,1,0)
\ee
for a constant $\lambda$ depending on the torsion type of the underlying structure,
\be
\lambda_\mrm{I}=2\,a^2,\quad\lambda_\mrm{II}=-4\,(a^2+b^2-c^2).
\ee
\begin{exa}
A cocalibrated $\G_2$-manifold of class $\mc{C}[\su(3)]$ with torsion type
$\mrm{I}$, for example, is homothetic to an $\eta$-Einstein-Sasakian manifold
whose Riemannian Ricci tensor is given by $\Ric^g=10\cdot g-4\cdot e_7\ox e_7$
(see \cite{Fri06}). A complete, simply connected, cocalibrated $\G_2$-manifold
of the same class but with torsion type $\mrm{II}$ is isometric \cite{Fri06} to
the product of a six-dimensional strictly nearly K\"ahler manifold with $\R$.
\end{exa}
After defining the $3$-form
\be
D_1:=(e_{127}+e_{347}+e_{567})-\omega^3=\diag(4,-4,0,0,0,0,0,0)
\ee
we can state the following result on necessary conditions for $\nabla^0\Psi=0$.
\begin{thm}
Let $(M^7,g,\omega^3)$ be a cocalibrated $\G_2$-manifold of class $\mc{C}[\su(3)]$
or $\mc{C}[\so(3)]$ with flux form $F=A\cdot\ast\omega^3$. Assuming $A\neq0$ there
exist a $\nabla^0$-parallel spinor field $\Psi_0$ if and only if the following is
satisfied:
\begin{enumerate}
  \item The spinor field $\Psi_0$ is $\nabla^c$-parallel and satisfies
        $D_1\cdot\Psi=4\cdot\Psi$.
  \item The parameter $B$ equals $1$.
\end{enumerate}
If $A=0$ instead, $\nabla^0\Psi_0=0$ implies either
\begin{enumerate}
  \item $\Ric^c\neq0$, $B=1$, $\Psi_0=\Psi_0^++\Psi_0^-$ and 
	      $D_1\cdot\Psi_0^\pm=\pm 4\cdot\Psi_0^\pm$ or
  \item $\Ric^c=0$ and $B=\pm1$.
\end{enumerate}
\end{thm}
\begin{rmk}
The condition $\Ric^c=0$ can only be realized for structures of class
$\mc{C}[\so(3)]$ with torsion type $\mrm{II}$.
\end{rmk}
A direct computation yields that two $\nabla^c$-parallel spinor
fields exist,
\be
\Psi_+ := [1,0,0,0,0,0,0,0]^T,\quad\Psi_- := [0,1,0,0,0,0,0,0]^T.
\ee
\begin{prop}
Any cocalibrated $\G_2$-structures of class $\mc{C}[\su(3)]$ or
$\mc{C}[\so(3)]$ admits two $\nabla^c$-parallel spinor fields
$\Psi_\pm$ such that $D_1\cdot\Psi_\pm=\pm 4\cdot\Psi_\pm$.
\end{prop}
Let us construct solutions to the non-special type of Killing spinor equation
based on these $\nabla^c$-parallel spinor fields but putting a more general
assumption on the $\nabla^c$-parallel flux form,
\be
F=A_1\cdot F_1+A_2\cdot (F_2+F_3),\quad A_1,A_2\in \R.
\ee
Here the $4$-forms $F_i$ are defined by
\be
F_1:=-e_{2467}+e_{2357}+e_{1457}+e_{1367},\quad F_2:=e_{1256}+e_{3456},\quad F_3:=e_{1234}.
\ee
\begin{thm}
Let $(M^7,g,\omega^3)$ be a cocalibrated $\G_2$-manifold of class $\mc{C}[\su(3)]$
or $\mc{C}[\so(3)]$ with flux form $F=A_1\cdot F_1+A_2\cdot (F_2+F_3)$. The equation
\bea
(\ast)&& s\,(X\haken T^c)\cdot\Psi_0 + p\,(X\haken F)\cdot\Psi_0 + q(X\hut F)\cdot\Psi_0=0
\eea
holds for the ansatz $\Psi_0=p_+\cdot\Psi_++p_-\cdot\Psi_-$ if and only if it
is solved by $p_+\cdot\Psi_+$ and $p_-\cdot\Psi_-$ separately, and if $a\cdot
s=b\cdot s=0$ in case of torsion type $\mrm{II}$. The spinor field $\Psi_\pm$
is a solution to $(\ast)$ if and only if the system
\be
(p+q)\,(\pm A_2-A_1)=\alpha^\pm\,s,\quad 4\,p\,A_1=\pm3\,q\,A_2-\beta^\pm\,s
\ee
holds for the parameters $\alpha^\pm$ and $\beta^\pm$, which depend on the torsion
type,
\be
\alpha^\pm_\mrm{I}=\pm a,\quad\alpha^\pm_\mrm{II}=-c,
\quad\beta^\pm_\mrm{I}=\pm3\,a,\quad\beta^\pm_\mrm{II}=0 .
\ee 
\end{thm}
\subsection{Cocalibrated \texorpdfstring{$\G_2$}{G2}-structures of class
\texorpdfstring{$\mc{C}[\su(2)]\cup\mc{C}[\un(2)]$}{C[su(2)], C[u(2)]}.}
In an adapted frame the three admissible types of torsion forms read
\bea
T^c_\mrm{I} &=& a\cdot( e_{127}+e_{347} ) + b\cdot e_{567} ,\\
T^c_\mrm{II}&=& a\cdot( e_{135}-e_{146}-e_{236}-e_{245})+b\cdot( e_{127}+e_{347}-2\,e_{567} )
,\quad a\neq0,\\
T^c_\mrm{III}&=& a\cdot( e_{135}-e_{245} )
\eea
for real parameters $a,b$. For structures of class $\mc{C}[\su(2)]$ with torsion type
$\mrm{I}$ or $\mrm{II}$ the relation $a=0$ or $a=-b$ is satisfied, respectively. There
exists no structure of class $\mc{C}[\un(2)]$ with torsion type $\mrm{III}$. The Ricci
tensor $\Ric^c$ of the characteristic connection is diagonal,
\be
\Ric^c=\diag(\lambda,\lambda,\lambda,\lambda,\kappa,\kappa,0).
\ee
Here the numbers $\lambda$ and $\kappa$ depend on the torsion type,
\be
\lambda_\mrm{I}=a^2+a\,b,\quad
\lambda_\mrm{II}=4\,a^2-b^2,\quad
\lambda_\mrm{III}=a^2,
\ee
\be
\kappa_\mrm{I}=2\,a\,b,\quad
\kappa_\mrm{II}=4\,a^2-4\,b^2,\quad
\kappa_\mrm{III}=0.
\ee
\begin{exa}
Any complete, simply connected, cocalibrated $\G_2$-manifold of class $\mc{C}[\su(2)]$
with torsion type $\mrm{I}$ -- for instance -- splits \cite{Fri06} into the product
$M^7=Y^4\x S^3$ of the sphere $S^3$ with a four-dimensional, complete, simply connected,
Ricci-flat, anti-selfdual manifold $Y^4$.
\end{exa}
For the formulation of classification results on $\nabla^0\Psi=0$ we split the
tangent bundle into
\be
TM^7=E_1\op E_2
\ee
where $E_2$ is spanned by $\{e_5,e_6,e_7\}$, and define the $3$-form
\bea
D_2&:=&-(e_{135}-e_{146}-e_{236}-e_{245})-\frac{1}{2}\cdot(e_{127}+e_{347}-2\,e_{567})\\
 & =&\diag(4,-4,-2,-2,1,1,1,1).
\eea
\begin{thm}
Let $(M^7,g,\omega^3)$ be a cocalibrated $\G_2$-manifold of class $\mc{C}[\su(2)]$
or $\mc{C}[\un(2)]$ with flux form $F=A\cdot\ast\omega^3$. Assuming $A\neq0$ there
exists a $\nabla^0$-parallel spinor field $\Psi_0$ if and only if the following
conditions are satisfied:
\begin{enumerate}
  \item The spinor field $\Psi_0$ is parallel with respect to $\nabla^c$ and fixed
        by $D_2\cdot\Psi_0=4\cdot\Psi_0$.
  \item $T$ coincides with the torsion form.
\end{enumerate}
If $A=0$ instead, $\nabla^0\Psi_0=0$ implies either
\begin{enumerate}
  \item $\Ric^c|_{E_1}\neq0$, $\Ric^c|_{E_2}\neq0$, $B=1$,
        $\Psi_0=\Psi_0^++\Psi_0^-$ and $D_2\cdot\Psi_0^\pm=\pm4\cdot\Psi_0^\pm$ or
  \item $\Ric^c|_{E_1}\neq0$, $\Ric^c|_{E_2}=0$, $B=1$,
        $\Psi_0=\Psi_0^++\Psi_0^-+\Psi_0^2$ and $D_2\cdot\Psi_0^2=-2\cdot\Psi_0^2$ or
  \item $\Ric^c=0$ and $B=\pm1$.
\end{enumerate}
\end{thm}
\begin{rmk}
The condition $\Ric^c=0$ can only be realized on certain structures of class
$\mc{C}[\un(2)]$ with torsion type $\mrm{I}$.
\end{rmk}
Again, the solution of $\nabla^0\Psi=0$ reduces to $\nabla^c\Psi=0$. A direct
computation leads to the existence of two spinor fields solving the latter
equation, $\Psi_+$ and $\Psi_-$ (see \autoref{ssec:1}). There exist another
two $\nabla^c$-parallel spinor fields for structures of class $\mc{C}[\su(2)]$,
\be
\Psi_1 := [0,0,1,0,0,0,0,0]^T,\quad\Psi_2 := [0,0,0,1,0,0,0,0]^T.
\ee
\begin{prop}
Any cocalibrated $\G_2$-structures of class $\mc{C}[\su(2)]$ or $\mc{C}[\un(2)]$
admits two $\nabla^c$-parallel spinor fields $\Psi_\pm$ such that $D_2\cdot
\Psi_\pm=\pm4\cdot\Psi_\pm$. There exist another two $\nabla^c$-parallel spinor
fields $\Psi_1,\Psi_2$ for structures of class $\mc{C}[\su(2)]$ satisfying
$D_2\cdot\Psi=-2\cdot\Psi$.
\end{prop}
Recall the definition of the $4$-forms $F_1$, $F_2$ and $F_3$ of the last subsection.
We proceed by constructing solutions to the non-special type of Killing spinor equation
based on the $\nabla^c$-parallel spinor fields $\Psi_\pm$. 
\begin{thm}
Let $(M^7,g,\omega^3)$ be a cocalibrated $\G_2$-manifold of class $\mc{C}[\su(2)]$ or
$\mc{C}[\un(2)]$ with flux form $F=A_1\cdot F_1+A_2\cdot F_2+A_3\cdot F_3$. The equation
\bea
(\ast)&& s\,(X\haken T^c)\cdot\Psi_0 + p\,(X\haken F)\cdot\Psi_0 + q(X\hut F)\cdot\Psi_0=0
\eea
holds for the ansatz $\Psi_0=p_+\cdot\Psi_++p_-\cdot\Psi_-$ if and only if it is
solved by $p_+\cdot\Psi_+$ and $p_-\cdot\Psi_-$ separately, and if $s=0$ in case
of torsion type $\mrm{III}$. $\Psi_\pm$ is a solution to $(\ast)$ if and only if
\be
(p+q)\,(\pm A_3-A_1)=\alpha^\pm\,s,\quad (p+q)\,(\pm A_2-A_1)=\beta^\pm\,s,\quad
4\,p\,A_1=\pm q\,(2\,A_2+A_3)-\gamma^\pm\,s.
\ee
The parameters $\alpha^\pm$, $\beta^\pm$ and $\gamma^\pm$ depend on the type of
torsion form,
\be
\alpha^\pm_\mrm{I}=\pm b,\,\,\alpha^\pm_\mrm{II}=-(\pm2\,b+a),
\,\,\beta^\pm_\mrm{I}=\pm a,\,\,\beta^\pm_\mrm{II}=\pm b-a,
\,\,\gamma^\pm_\mrm{I}=\pm(2\,a+b)\,s,\,\,\gamma^\pm_\mrm{II}=0.
\ee
\end{thm}\noindent
We then proceed with construction results based on $\Psi_1$ and $\Psi_2$.
\begin{thm}
Let $(M^7,g,\omega^3)$ be a cocalibrated $\G_2$-manifold of class $\mc{C}[\su(2)]$
with flux form $F=A_1\cdot F_1+A_2\cdot F_2+A_3\cdot F_3$. The equation
\bea
(\ast)&& s\,(X\haken T^c)\cdot\Psi_0 + p\,(X\haken F)\cdot\Psi_0 + q(X\hut F)\cdot\Psi_0=0
\eea
holds for the ansatz $\Psi_0=p_1\cdot\Psi_1+p_2\cdot\Psi_2$ if and only if it is
solved by $p_1\cdot\Psi_1$ and $p_2\cdot\Psi_2$ separately. The spinor field $\Psi_i$
is a solution to $(\ast)$ if and only if
\be
2\,(p+q)\,A_1=\alpha^i\,s,\quad 2\,(p+q)\,A_2=\beta^i\,s,\quad
2\,(p+q)\,A_3=\gamma^i\,s,\quad 2\,p\,A_2+q\,A_3=\delta^i\,s.
\ee
The numbers $\alpha^i$, $\beta^i$, $\gamma^i$ and $\delta^i$ depend on the torsion type,
\be
\alpha^i_\mrm{I}=0,\quad\alpha^i_\mrm{II}=2\,a,\quad\alpha^i_\mrm{III}=a,
\quad\beta^i_\mrm{I}=2\,a,\quad\beta^i_\mrm{II}=2\,b,\quad\beta^i_\mrm{III}=(-1)^{i-1}a,
\ee
\be
\gamma^i_\mrm{I}=2\,b,\quad\gamma^i_\mrm{II}=-4\,b,\quad\gamma^i_\mrm{III}=(-1)^{i-1}a,
\quad\delta^i_\mrm{I}=b,\quad\delta^i_\mrm{II}=-2\,b,\quad\delta^i_\mrm{III}=(-1)^{i}a.
\ee
\end{thm}
\begin{rmk}
This construction shows that for cocalibrated $\G_2$-structures of class $\mc{C}[\su(2)]$
with torsion type $\mrm{I}$ or $\mrm{II}$ there exist at most four linearly independent
$\nabla$-parallel spinor fields ($F\neq0$), with type $\mrm{III}$ at most three.
\end{rmk}
\begin{exa}
Take a structure of class $\mc{C}[\su(2)]$ with torsion type $\mrm{II}$. The four
spinor fields $\Psi_+$, $\Psi_-$, $\Psi_1$ and $\Psi_2$ solve $(\ast)$ if and only
if
\be
p=0,\quad q\neq0,\quad q\,A_1=a\,s,\quad q\,A_2=b\,s,\quad A_3=-2\,A_2.
\ee
\end{exa}
\subsection{Cocalibrated \texorpdfstring{$\G_2$}{G2}-structures of class
\texorpdfstring{$\mc{C}[\su_c(2)\,\mrm{rel.}]$}{C[suc(2) rel.]}.}
Using an adapted frame the torsion form is given by
\be
T^c=a\cdot\omega^3+b\cdot e_{567}, \quad b\neq0
\ee
for real parameters $a$, $b$. The characteristic Ricci tensor has diagonal
form,
\be
\Ric^c=\diag(\lambda,\lambda,\lambda,\lambda,\kappa,\kappa,\kappa),
\ee
and the numbers $\lambda$ and $\kappa$ are
\be
\lambda=12\,a^2+3\,a\,b,\quad\kappa=12\,a^2+4\,a\,b.
\ee
\begin{exa}
Any complete, simply connected, cocalibrated $\G_2$-manifold of class $\mc{C}[\su_c(2)]$ or
$\mc{C}[\R\op\su_c(2)]$ is a naturally reductive homogeneous space (see \cite{Fri06}).
\end{exa}
We recall the splitting of the tangent bundle $TM^7=E_1\op E_2$ introduced in the last
subsection, define the $3$-form
\be
D_3:=e_{567}-\omega^3=\diag(6,-2,-2,-2,0,0,0,0)
\ee
and state the classification results regarding $\nabla^0\Psi=0$.
\begin{thm}
Let $(M^7,g,\omega^3)$ be a cocalibrated $\G_2$-manifold of class $\mc{C}[\su_c(2)\,
\mrm{rel.}]$ with flux form $F=A\cdot\ast\omega^3$. If there exists a $\nabla^0$-parallel
spinor field $\Psi_0$ for $A\neq0$, one of the following holds:
\begin{enumerate}
  \item The spinor field is $\nabla^c$-parallel and satisfies $D_3\cdot\Psi=6\cdot\Psi$.
        The parameter $B$ equals $1$.
  \item The spinor field is fixed by $D_3\cdot\Psi=-2\cdot\Psi$. The system parameters
        satisfy $A=-2\,a$, $B=-7$, $b=3\,a$.
\end{enumerate}
If $A=0$ instead, $\nabla^0\Psi_0=0$ implies one of the following:
\begin{enumerate}
  \item $\Ric^c|_{E_1}\neq0$, $\Ric^c|_{E_2}\neq0$, $B=1$
        and $D_3\cdot\Psi_0=6\cdot\Psi_0$.
  \item $\Ric^c|_{E_1}\neq0$, $\Ric^c|_{E_2}=0$, $B=1$,
        $\Psi_0=\Psi_0^1+\Psi_0^3$ and $D_3\cdot\Psi_0^3=-2\cdot\Psi_0^3$.
  \item $\Ric^c=0$ and $B=\pm1$.
\end{enumerate}
\end{thm}
A direct computation leads to the existence of a $\nabla^c$-parallel spinor field,
namely $\Psi_+$ (see \autoref{ssec:1}).
\begin{prop}
Cocalibrated $\G_2$-structures of class $\mc{C}[\su_c(2)\,\mrm{rel.}]$ admit
a $\nabla^c$-parallel spinor field $\Psi_+$ satisfying $D_3\cdot\Psi=6\cdot\Psi$.
\end{prop}
We conclude with construction results regarding the non-special type of Killing spinor
equation.
\begin{thm}
Let $(M^7,g,\omega^3)$ be a cocalibrated $\G_2$-manifold of class $\mc{C}[\su_c(2)\,
\mrm{rel.}]$ with flux form $F=A_1\cdot (F_1+F_2)+A_2\cdot F_3$. The equation
\be
s\,(X\haken T^c)\cdot\Psi_+ + p\,(X\haken F)\cdot\Psi_+ + q(X\hut F)\cdot\Psi_+=0
\ee
holds if and only if
\be
(p+q)\,(A_1-A_2)=-b\,s,\quad 3\,(p-q)\,A_1+p\,A_2=-3\,a\,s.
\ee
\end{thm}
%
%
%
%
\section{Conclusions}\label{sec:8}\noindent
Let us go back to the full system $(\aleph)$ \& $(\aleph\aleph)$. If we define
\be
\Ric^T:=\Ric^g_{ij} - \frac{1}{4}\,T_{imn}T_{jmn},
\ee
the relations
\be
\delta T=0,\quad F\cdot\Psi=\kappa\cdot\Psi,\quad
\mrm{div}^c(\Ric^T)=\mrm{div}(\Ric^T)=0
\ee
hold in all cases discussed here. If we replace -- as suggested
in \cite{AFNP05} -- the equation $\Ric^T=0$ in $(\aleph)$ by $\mrm{div}(\Ric^T)
=0$, then every constructed solution satisfies the new system, provided
furthermore that the spinor field $\Psi$ is an eigenspinor of $T$.

\autoref{tab:1} contains a summary of our results. For instance, cocalibrated
$\G_2$-structures of class $\mc{C}[\su(2)]$ with torsion type I will admit at
most four distinct (i.e.~linearly independent) spinor fields parallel with respect
to a certain family of spinorial covariant derivatives with non-vanishing $4$-form.
Should we further consider $T=T^c$, there exists a family $\bar{\nabla}$ rendering
three distinct spinor fields parallel. At the same time $\Ric^{T^c}=\Ric^c=0$ is
fulfilled. There exist other structures of the same type and a family $\tilde{\nabla}$
with now four distinct spinor fields such that $\tilde{\nabla}\tilde{\Psi}=0$,
$\Ric^T=0$ but $T\neq T^c$. The `check' mark indicates that all these spinor fields
are eigenspinors relative to the $3$-form $T$. Finally, there exist four distinct
$\nabla^c$-parallel spinor fields.

To conclude, a few comments on possible generalizations of our spinorial
field equations. Solving the following Killing spinor equation is the main concern
of supergravity models in type II string theory (see \cite{GLW05}):
\be
\nabla^g_X\Psi+\frac{1}{4}\,(X\haken T)\cdot\Psi+\sum_i p_i\,(X\haken F^i)\cdot\Psi
+\sum_i q_i\,(X\hut F^i)\cdot\Psi=0.
\ee
The differential forms $F^i$ are of degree $2\,i$ (type IIa) or of
degree $2\,i+1$ (type IIb). Due to the complexity of the algebraic systems
the approach of this paper is unlikely to be suitable for this general
kind of equation. However, $\nabla^c$-parallel spinor fields may represent
natural candidates to begin with when constructing solutions. If we start
from one of the structures considered in this work and set $T=B\cdot T^c$,
then the above equation will read
\be
\frac{B-1}{4}\,(X\haken T^c)\cdot\Psi_0+\sum_i p_i\,(X\haken F^i)\cdot\Psi_0
+\sum_i q_i\,(X\hut F^i)\cdot\Psi_0=0
\ee
for a $\nabla^c$-parallel spinor field $\Psi_0$. The last column of
table \ref{tab:1} tells us how many such spinor fields exist. We
conjecture that this purely algebraic equation could be solved with
an appropriate ansatz for the differential forms $F^i$.
\renewcommand{\arraystretch}{1.3}
\begin{table}[t]
  \caption{\label{tab:1} Existence of solutions to $\nabla\Psi=0$.
           $\mc{N}$ denotes the maximum number of constructed, linearly
           independent spinor fields which are parallel with respect to a
           certain family of spinorial covariant derivatives $\nabla$.
           The superscript $c$ refers to the characteristic connection
           $\nabla^c$. In the second last column we determine whether
           all constructed solutions are eigenspinors of the differential
           form $T$.}
  {\small
  \begin{center}
  \begin{tabular}{|c|ccc||c|c|c|c||c|}
   \hline
   Dim. & \multicolumn{3}{c||}{Structure} & \multicolumn{4}{c||}{$F\neq0$}
        & $F=0$ \\ \cline{5-9}
        &&&& $\mc{N}$ & \multicolumn{2}{c|}{$\mc{N}(\Ric^{T}=0)$}
        & $T\cdot\Psi=\lambda\cdot\Psi$ & $\mc{N}^c$ \\ \cline{6-7}
        &&&&& $T\neq T^c$ & $T=T^c$ & & \\
   \hline \hline
   $n=5$ & \multicolumn{3}{c||}{$\alpha$-Sasakian structure}
         &$1$&$1$&$1$&\ding{52}&$2$ \\
   \hline \hline
   $n=6$ & \multicolumn{1}{c|}{almost}   & \multicolumn{2}{c||}{$\SU(3)$}
         &$2$&$1$&--&\ding{52}&$2$ \\ \cline{3-9}
         & \multicolumn{1}{c|}{Hermitian}   & \multicolumn{2}{c||}{$\SO(3)$}
         &$2$&$1$&$2$&\ding{52}&$2$ \\ \cline{3-9}
         & \multicolumn{1}{c|}{structure}  & \multicolumn{2}{c||}{$\SU(2)$}
         &$2$&--&--&\ding{56}&$4$ \\ \cline{3-9}
         & \multicolumn{1}{c|}{of}    & \multicolumn{2}{c||}{$\U(2)_0$}
         &$2$&--&--&\ding{56}&$4$  \\ \cline{3-9}
         & \multicolumn{1}{c|}{class} & \multicolumn{2}{c||}{$\U(2)_1$}
         &\multicolumn{5}{l|}{no solutions}\\ \cline{3-9}
         & \multicolumn{1}{c|}{$\mc{C}[\G]$} & \multicolumn{2}{c||}{$\U(2)_{-1}$}
         &$2$&$1$&--&\ding{52}&$2$ \\
   \hline \hline
   $n=7$ & \multicolumn{3}{c||}{nearly parallel $\G_2$-structure}
         &$2$&$2$&--&\ding{52}&$1$ \\ \cline{2-9}
         & \multicolumn{1}{c|}{cocalibrated} & \multicolumn{1}{c|}{$\su(3)$} & I
         &$2$&--&--&\ding{52}&$2$ \\ \cline{4-9}
         & \multicolumn{1}{c|}{$\G_2$-structure} & \multicolumn{1}{c|}{} & II
         &$2$&$2$&--&\ding{52}&$2$ \\ \cline{3-9}
         & \multicolumn{1}{c|}{of}    & \multicolumn{1}{c|}{$\so(3)$} & I
         &$2$&--&--&\ding{52}&$2$ \\ \cline{4-9}
         & \multicolumn{1}{c|}{class} & \multicolumn{1}{c|}{} & II
         &$2$&$2$&$1$&\ding{52}&$2$ \\ \cline{3-9}
         & \multicolumn{1}{c|}{$\mc{C}[\g]$} & \multicolumn{1}{c|}{$\su(2)$} & I
         &$4$&$4$&$3$&\ding{52}&$4$ \\ \cline{4-9}
         & \multicolumn{1}{c|}{} & \multicolumn{1}{c|}{}& II
         &$4$&--&--&\ding{52}&$4$ \\ \cline{4-9}
         & \multicolumn{1}{c|}{} & \multicolumn{1}{c|}{}& III
         &$3$&--&--&\ding{52}&$4$ \\ \cline{3-9}
         & \multicolumn{1}{c|}{} & \multicolumn{1}{c|}{$\un(2)$}& I
         &$2$&$2$&$2$&\ding{52}&$2$ \\ \cline{4-9}
         & \multicolumn{1}{c|}{} & \multicolumn{1}{c|}{}& II
         &$2$&$2$&--&\ding{52}&$2$ \\ \cline{3-9}
         & \multicolumn{1}{c|}{} & \multicolumn{2}{c||}{$\su_c(2)\,\mrm{rel.}$ }
         &$1$&$1$&$1$&\ding{52}&$1$ \\ \cline{2-9}
         & \multicolumn{3}{c||}{$\alpha$-Sasakian structure}
         &$2$&$2$&$2$&\ding{52}&$2$ \\ \hline
  \end{tabular}
  \vspace{0.5cm}
  \end{center}}
\end{table}
\renewcommand{\arraystretch}{1.0}
\begin{thank}
We wish to thank Thomas Friedrich and Nils Schoemann for discussions
and both the SPP 1154:~\emph{Global Differential Geometry} and the
SFB 647:~\emph{Space--Time--Matter} for financial support.
\end{thank}
%
%
%
%

%
%
%

\begin{bibdiv}
%
\begin{biblist}
%
\bib{AF03}{article}{
  author={Agricola, I.},
  author={Friedrich, T.},
  title={Killing spinors in supergravity with {$4$}-fluxes},
  journal={Class. Quant. Grav.},
  volume={20},
  date={2003},
  pages={4707--4717},
}
%
\bib{AFNP05}{article}{
  author={Agricola, I.},
  author={Friedrich, T.},
  author={Nagy, P.-A.},
  author={Puhle, C.},
  title={On the Ricci tensor in the common sector of type II string theory},
  journal={Class. Quant. Grav.},
  volume={22},
  date={2005},
  pages={2569--2577},
}
%
\bib{AFS05}{article}{
  author={Alexandrov, B.},
  author={Friedrich, T.},
  author={Schoemann, N.},
  title={Almost Hermitian {$6$}-Manifolds Revisited},
  journal={Journ. Geom. Phys.},
  volume={53},
  date={2005},
  pages={1--30},
}
%
\bib{AG86}{article}{
  author={Aleksiev, V.},
  author={Ganchev, G.},
  title={On the classification of the almost contact metric manifolds},
  journal={in: Proc. 15th Spring Conf., Sofia, Bulgaria, UBM},
  date={1986},
  pages={155--161},
}
%
\bib{Ali01}{article}{
  author={Ali, T.},
  title={{$\M $}-theory on seven manifolds with {$\G$}-fluxes},
  pages={hep-th/0111220},
}
%
\bib{BFGK91}{book}{
  author={Baum, H.},
  author={Friedrich, T.},
  author={Grunewald, R.},
  author={Kath, I.},
  title={Twistors and Killing
         spinors on Riemannian manifolds},
  series={Teubner-Texte zur Mathematik},
  publisher={Teubner-Verlag},
  address={Leipzig},
  volume={124},
  date={1991},
}
%
\bib{BDS01}{article}{
  author={Bilal, A.},
  author={Derendinger, J.-P.},
  author={Sfetsos, K.},
  title={Weak {$\G _2$}-holonomy from self-duality, flux and supersymmetry},
  journal={Nucl. Phys. B},
  volume={628},
  date={2002},
  pages={112--132},
}
%
\bib{Bla76}{book}{
  author={Blair, D. E.},
  title={Contact manifolds in Riemannian geometry},
  series={Lecture Notes in Mathematics},
  volume={509},
  publisher={Springer},
  date={1976},
}
%
\bib{BGN00}{article}{
  author={Boyer, C. P.},
  author={Galicki, K.},
  author={Nakamaye, M.},
  title={On the geometry of Sasakian-Einstein {$5$}-manifolds},
  journal={Math. Ann.},
  volume={325},
  date={2003},
  pages={485--524},
}
%
\bib{BGM04}{article}{
  author={Boyer, C. P.},
  author={Galicki, K.},
  author={Matzeu, P.},
  title={On eta-Einstein sasakian geometry},
  journal={Comm. Math. Phys.},
  volume={262},
  date={2006},
  pages={177--208},
}
%
\bib{Bry87}{article}{
  author={Bryant, R. L.},
  title={Metrics with exceptional holonomy},
  journal={Ann. of Math. (2)},
  volume={126},
  date={1987},
  pages={525--576},
} 
%
\bib{But05}{article}{
  author={Butruille, J. B.},
  title={Classification des variétés approximativement
         k\"ahleriennes homogènes},
  journal={Ann. Global Anal. Geom.},
  volume={27},
  date={2005},
  pages={201--225},
}
%
\bib{CG90}{article}{
  author={Chinea, D.},
  author={Gonzales, G.},
  title={A classification of almost contact metric
         manifolds},
  journal={Ann. di mat. Pura Appl.},
  volume={156},
  date={1990},
  pages={15--36},
}
%
\bib{Duf02}{article}{
  author={Duff, M. J.},
  title={{$\M $}-theory on manifolds of {$\G _2$}-holonomy: the first twenty years},
  pages={hep-th/0201062},
}
%
\bib{Dyn57}{article}{
  author={Dynkin, E. B.},
  title={The maximal subgroups of the classical groups},
  journal={Am. Math. Soc., Transl., II},
  volume={6},
  date={1957},
  pages={245--378}
}
%
\bib{FG82}{article}{
  author={Fernandez, M.},
  author={Gray, A.},
  title={Riemannian manifolds with structure
         group {$\G_2$}},
  journal={Ann. Mat. Pura Appl.},
  volume={132},
  date={1982},
  pages={19--45},
}
%
\bib{Fri00}{book}{
  author={Friedrich, T.},
  title={Dirac operators in Riemannian geometry},
  series={Graduate Studies in Mathematics},
  volume={25},
  publisher={AMS},
  date={2000},
  address={Providence},
}
%
\bib{Fri06}{article}{
  author={Friedrich, T.},
  title={{$\G_2$}-manifolds with parallel characteristic torsion},
  journal={Diff. Geom. Appl.},
  volume={25},
  date={2007},
  pages={632--648},
}
%
\bib{FI03a}{article}{
  author={Friedrich, T.},
  author={Ivanov, S.},
  title={Almost contact manifolds, connections with
         torsion and parallel spinors},
  journal={J. Reine Angew. Math.},
  volume={559},
  date={2003},
  pages={217--236},
}
%
\bib{FI02}{article}{
  author={Friedrich, T.},
  author={Ivanov, S.},
  title={Parallel spinors and connections with skew-symmetric torsion in string theory},
  journal={Asian Journ. Math},
  volume={6},
  date={2002},
  pages={303--336},
}
%
\bib{FK89}{article}{
  author={Friedrich, T.},
  author={Kath, I.},
  title={Einstein manifolds of dimension five
         with small eigenvalues of the Dirac operator},
  journal={Journ. Diff. Geom.},
  volume={19},
  date={1989},
  pages={263--279},
}
%
\bib{FK90}{article}{
  author={Friedrich, T.},
  author={Kath, I.},
  title={Compact seven-dimensional Riemannian manifolds
         with Killing spinors},
  journal={Comm. Math. Phys.},
  volume={133},
  date={1990},
  pages={543--561},
}
%
\bib{FKMS97}{article}{
  author={Friedrich, T.},
  author={Kath, I.},
  author={Moroianu, A.},
  author={Semmelmann, U.},
  title={On nearly parallel {$\G_2$}-structures},
  journal={Journ. Geom. Phys.},
  volume={23},
  date={1997},
  pages={256--286},
}
%
\bib{FK00}{article}{
  author={Friedrich, T.},
  author={Kim, E. C.},
  title={The Einstein-Dirac equation on Riemannian
         spin manifolds},
  journal={Journ. Geom. Phys.},
  volume={33},
  date={2000},
  pages={128--172},
}
%
\bib{GLW05}{article}{
  author={Gra\~na, M.},
  author={Louis, J.},
  author={Waldram, D.},
  title={Hitchin Functionals in {$N = 2$} Supergravity},
  journal={JHEP},
  volume={0601},
  date={2006},
  pages={008},
}
%
\bib{Gra76}{article}{
  author={Gray, A.},
  title={The structure of nearly K\"ahler manifolds},
  journal={Math. Ann.},
  volume={223},
  date={1976},
  pages={233--248},
}
%
\bib{GH80}{article}{
  author={Gray, A.},
  author={Hervella, L.},
  title={The sixteen classes of almost Hermitian manifolds
         and their linear invariants},
  journal={Ann. Mat. Pura Appl.},
  volume={123},
  date={1980},
  pages={35--58},
}
%
\bib{Puh07}{thesis}{
  author={Puhle, C.},
  title={Spinorial field equations in supergravity with fluxes},
  eprint={http://www.math.hu-berlin.de/~puhle/publications.html},         
  organization={Hum\-boldt University of Berlin},
  date={2007},
  type={Ph.D. Thesis},
}
%
\bib{Sch06}{thesis}{
  author={Schoemann, N.},
  title={Almost Hermitian structures with parallel torsion},
  organization={Hum\-boldt University of Berlin},
  date={2006},
  type={Ph.D. Thesis},
}
%
\bib{Sch07}{article}{
  author={Schoemann, N.},
  title={Almost hermitian structures with parallel torsion},
  journal={J. Geom. Phys.},
  volume={57},
  date={2007},
  pages={2187--2212}, 
}
%
\bib{Str86}{article}{
  author={Strominger, A.},
  title={Superstrings with torsion},
  journal={Nucl. Phys. B},
  volume={274},
  date={1986},
  pages={253--284},
}
%
\bib{WNW85}{article}{
  author={Witt, B. de},
  author={Nicolai, H.},
  author={Warner, N. P.},
  title={The embedding of gauged {$n=8$} supergravity into {$d=11$} supergravity},
  journal={Nucl. Phys. B},
  volume={255},
  date={1985},
  pages={29},
}
%
\end{biblist}
%
\end{bibdiv}
\end{document}